\documentclass[aop,seceqn,nameyear,MSNbibl,dvips]{arximspdf}
\usepackage{mathrsfs}

\doi{10.1214/10-AOP629}
\volume{40}
\issue{2}
\pubyear{2012}
\firstpage{813}
\lastpage{857}

\makeatletter

\newtheorem{theorem}{Theorem}[section]
\newtheorem{lemma}[theorem]{Lemma}
\newtheorem{corol}[theorem]{Corollary}
\newtheorem{prop}[theorem]{Proposition}

\newcommand{\mcr}{\mathscr}
\newcommand{\mbb}{\mathbb}
\newcommand{\mbf}{\mathbf}

\newcommand{\const}{\mathrm{const}}

\makeatother

\begin{document}
\begin{frontmatter}

\title{Stochastic equations, flows and~measure-valued~processes}
\runtitle{Stochastic equations, flows and processes}

\begin{aug}
\author[A]{\fnms{Donald A.} \snm{Dawson}\thanksref{t1}\ead[label=e1]{ddawson@math.carleton.ca}\ead[label=u1,url]{http://lrsp.carleton.ca/directors/dawson/}}
and
\author[B]{\fnms{Zenghu} \snm{Li}\corref{}\thanksref{t2}\ead[label=e2]{lizh@bnu.edu.cn}\ead[label=u2,url]{http://math.bnu.edu.cn/\textasciitilde lizh/}}
\runauthor{D. A. Dawson and Z. Li}
\affiliation{Carleton University and Beijing Normal University}
\address[A]{School of Mathematics and Statistics\\
Carleton University\\
1125 Colonel By Drive\\
Ottawa, Ontario\\
Canada K1S 5B6 \\
\printead{e1}\\
\printead{u1}} 
\address[B]{School of Mathematical Sciences\\
Beijing Normal University\\
Beijing 100875\\
P. R. China\\
\printead{e2}\\
\printead{u2}}
\end{aug}

\thankstext{t1}{Supported by NSERC.}
\thankstext{t2}{Supported by NSFC and CJSP.}

\received{\smonth{1} \syear{2010}}
\revised{\smonth{9} \syear{2010}}

%
\begin{abstract}
We first prove some general results on pathwise uniqueness, comparison
property and existence of nonnegative strong solutions of stochastic
equations driven by white noises and Poisson random measures. The
results are then used to prove the strong existence of two classes of
stochastic flows associated with coalescents with multiple collisions,
that is, generalized Fleming--Viot flows and flows of continuous-state
branching processes with immigration. One of them unifies the different
treatments of three kinds of flows in Bertoin and Le Gall [\textit{Ann.
Inst. H. Poincar\'{e} Probab. Statist.} \textbf{41} (2005) 307--333]. Two
scaling limit theorems for the generalized Fleming--Viot flows are
proved, which lead to sub-critical branching immigration
superprocesses. From those theorems we derive easily a generalization
of the limit theorem for finite point motions of the flows in Bertoin
and Le Gall [\textit{Illinois J. Math.} \textbf{50} (2006) 147--181].
\end{abstract}

%

\end{frontmatter}

\section{Introduction}

A class of stochastic flows of bridges were introduced by Bertoin and Le
Gall (\citeyear{BeL03}) to study the coalescent processes with
multiple collisions of
\citet{Pit99} [see also \citet{Sag99}]. The law of such a coalescent
process is determined by a finite measure $\Lambda(dz)$ on $[0,1]$. The
Kingman coalescent corresponds to $\Lambda= \delta_0$ and the
Bolthausen--Sznitman coalescent corresponds to $\Lambda=$ Lebesgue
measure on $[0,1]$ [see \citet{BoS98} and \citet{Kin09}].
In fact,
\citet{BeL03} established a remarkable
connection between the coalescents with multiple collisions and the
stochastic flows of bridges. Based on this connection, they have
developed a theory of the coalescents and the flows in the series of
papers [see Bertoin and Le Gall (\citeyear{BeL03}, \citeyear{BeL05},
\citeyear{BeL06})]. We refer the reader
to \citet{LeR04}, \citet{MaX01} and \citet{Xia09} for
the study of stochastic flows of mappings and measures in abstract
settings.

Let $\{B_{s,t}\dvtx-\infty< s\le t<\infty\}$ be the stochastic flow of
bridges associated to a $\Lambda$-coalescent in the sense of
\citet{BeL03}.
A number of precise characterizations of the flow
$\{B_{-t,0}(v)\dvtx t\ge0, v\in[0,1]\}$ were given in \citet{BeL03}.
For any $t\ge0$, the function $v\mapsto B_{-t,0}(v)$ induces a
random probability measure $\rho_t(dv)$ on $[0,1]$. The process
$\{\rho_t\dvtx t\ge0\}$ was characterized in \citet{BeL03} as
the unique solution of a~martingale problem. In fact, this process is a
measure-valued dual to the $\Lambda$-coalescent process. It was also
pointed out in \citet{BeL03} that $\{\rho_t\dvtx t\ge0\}$ can
be regarded as a \textit{generalized Fleming--Viot process} [see also
Donnelly and Kurtz (\citeyear{DoK99b}, \citeyear{DoK99a})].

Let $\Lambda(dz)$ be a finite measure on $[0,1]$ such that
$\Lambda(\{0\}) = 0$, and let $\{M(ds,dz,du)\}$ be a Poisson random
measure on $(0,\infty)\times(0,1]^2$ with intensity
$z^{-2}\,ds\,\Lambda(dz)\,du$. It was proved in \citet{BeL05}
that there is weak solution flow $\{X_t(v)\dvtx t\ge0, v\in
[0,1]\}$ to the stochastic equation
%
%
\begin{equation}\label{1.1}
X_t(v) = v + \int_0^t\int_0^1\int_0^1 z\bigl[1_{\{u\le X_{s-}(v)\}} -
X_{s-}(v)\bigr] M(ds,dz,du).
\end{equation}
Moreover, \citet{BeL05} showed that for any $0\le r_1<
\cdots< r_p\le1$ the $p$-point motion $\{(B_{-t,0}(r_1), \ldots,
B_{-t,0}(r_p))\dvtx t\ge0\}$ is equivalent to $\{(X_t(r_1), \ldots,
X_t(r_p))\dvtx t\ge0\}$. Therefore, the solutions of (\ref{1.1}) give a
realization of the flow of bridges associated with the
$\Lambda$-coalescent process. A separate treatment for the Kingman
coalescent flow was also given in \citet{BeL05}. In
that case they showed the $p$-point motion $\{(B_{-t,0}(r_1),
\ldots, B_{-t,0}(r_p))\dvtx t\ge0\}$ is a diffusion process in
\[
D_p := \{x=(x_1,\ldots,x_p)\in\mbb{R}^p\dvtx0\le x_1\le\cdots\le
x_p\le1\}
\]
with generator $A_0$ defined by
%
%
\begin{equation}\label{1.2}
A_0f(x) = \frac{1}{2}\sum_{i,j=1}^p x_{i\land j}(1-x_{i\vee
j})\,\frac{\partial^2f}{\partial x_i\,\partial x_j}(x).
\end{equation}
Given a $\Lambda$-coalescent flow $\{B_{s,t}\dvtx-\infty< s\le
t<\infty\}$, we define the \textit{flow of inverses} by
\[
B_{s,t}^{-1}(v) = \inf\{u\in[0,1]\dvtx B_{s,t}(u)>v\}, \qquad v\in
[0,1),
\]
and $B_{s,t}^{-1}(1) = B_{s,t}^{-1}(1-)$. In the Kingman coalescent
case, it was proved in \citet{BeL05} that the $p$-point
motion $\{(B_{0,t}^{-1}(r_1), \ldots, B_{0,t}^{-1}(r_p))\dvtx\allowbreak t\ge0\}$
is a diffusion process in $D_p$ with generator $A_1$ given by
%
%
\begin{equation}\label{1.3}
A_1f(x) = A_0f(x) + \sum_{i=1}^p\biggl(\frac{1}{2}-x_i\biggr)\,
\frac{\partial f}{\partial x_i}(x),
\end{equation}
where $A_0$ is given by (\ref{1.2}). The analogous characterization for
the $\Lambda$-coa\-lescent flow with $\Lambda(\{0\}) = 0$ was also
provided in \citet{BeL05}. Those results give deep insights
into the structures of the stochastic flows associated with the
$\Lambda$-coalescents.

The asymptotic properties of $\Lambda$-coalescent flows were studied in
\citet{BeL06}. For each integer $k\ge1$ let
$\Lambda_k(dx)$ be a finite measure on $[0,1]$ with $\Lambda_k(\{
0\})
= 0$ and let $\{X_k(t,v)\dvtx t\ge0,v\in[0,1]\}$ be defined by (\ref{1.1})
from a Poisson random measure $\{M_k(ds,dz,du)\}$ on $(0,\infty)\times
(0,1]^2$ with intensity $z^{-2}\,ds\,\Lambda_k(dz)\,du$. Suppose that
$z^{-2}(z\land z^2) \Lambda_k(k^{-1}\,dz)$ converges weakly as $k\to
\infty$ to a finite measure on $(0,\infty)$ denoted by $z^{-2}(z\land
z^2) \Lambda(dz)$. By a limit theorem of \citet{BeL06} the
rescaled $p$-point motion $\{(kX_k(kt,r_1/k), \ldots,
kX_k(kt,r_p/k))\dvtx
t\ge0\}$ converges in distribution to that of the weak solution flow of
the stochastic equation
%
%
\begin{equation}\label{1.4}
Y_t(v) = v + \int_0^t\int_0^\infty\int_0^\infty x 1_{\{u\le
Y_{s-}(v)\}} \tilde{N}(ds,dx,du),
\end{equation}
where $\tilde{N}(ds,dx,du)$ is a compensated Poisson random measure on
$[0,\infty)\times(0,\infty)^2$ with intensity $z^{-2}\,ds\,\Lambda(dz)\,du$.
It was pointed out in \citet{BeL06} that the solution of
(\ref{1.4}) is a special critical \textit{continuous-state branching
process} (CB-process).

In this paper we study two classes of stochastic flows defined by
stochastic equations that generalize (\ref{1.1}) and (\ref{1.4}). We
shall first treat the generalization of (\ref{1.4}) since it involves
simpler structures. Suppose that $\sigma\ge0$ and~$b$ are constants,
$v\mapsto\gamma(v)$ is a nonnegative and nondecreasing continuous
function on $[0,\infty)$ and $(z\land z^2)m(dz)$ is a finite measures on
$(0,\infty)$. Let $\{W(ds,du)\}$ be a white noise on $(0,\infty)^2$ based
on the Lebesgue measure\vspace*{1pt} $ds\,du$. Let $\{N(ds,dz,du)\}$ be
a Poisson random
measure on $(0,\infty)^3$ with intensity $ds\,m(dz)\,du$. Let
$\{\tilde{N}(ds,dz,du)\}$ be the compensated measure of $\{N(ds,dz,du)\}$. We
shall see that for any $v\ge0$ there is a pathwise unique nonnegative
solution of the stochastic equation
%
%
\begin{eqnarray}\label{1.5}
Y_t(v) &=& v + \sigma\int_0^t\int_0^{Y_{s-}(v)} W(ds,du) +
\int_0^t [\gamma(v)-bY_{s-}(v)] \,ds \nonumber\\[-8pt]\\[-8pt]
&&{}
+ \int_0^t \int_0^\infty\int_0^{Y_{s-}(v)} z
\tilde{N}(ds,dz,du).\nonumber
\end{eqnarray}
It is not hard to show each solution $Y(v) = \{Y_t(v)\dvtx t\ge0\}$ is a
\textit{continuous-state branching process with immigration}
(CBI-process). Then it is natural to call the two-parameter process
$\{Y_t(v)\dvtx t\ge0, v\ge0\}$ a \textit{flow of CBI-processes}. We prove
that the flow has a version with the following properties:
\begin{longlist}
\item
for each $v\ge0$, $t\mapsto Y_t(v)$ is a c\`{a}dl\`{a}g process on
$[0,\infty)$ and solves~(\ref{1.5});

\item
for each $t\ge0$, $v\mapsto Y_t(v)$ is a nonnegative and
nondecreasing c\`{a}dl\`{a}g process on $[0,\infty)$.\vadjust{\goodbreak}
\end{longlist}
The proof of those properties is based on the observation that
$\{Y(v)\dvtx
v\ge0\}$ is a path-valued process with independent increments. For any
$t\ge0$, the random function $v\mapsto Y_t(v)$ induces a random Radon
measure $Y_t(dv)$ on $[0,\infty)$. We shall see that $\{Y_t\dvtx t\ge0\}
$ is
actually an \textit{immigration superprocess} in the sense of
\citet{Li10}
with trivial underlying spatial motion. One could replace the diffusion
term in (\ref{1.5}) by the stochastic integral $\sigma\int_0^t
\sqrt{Y_{s-}(v)} \,dW(s)$ using a one-dimensional Brownian motion
$\{W(t)\dvtx
t\ge0\}$ as in \citet{DaL06}. The resulted equation defines an
equivalent CBI-process for any fixed $v\ge0$, but it does not give an
equivalent flow.

To describe our generalization of (\ref{1.1}), let us assume that
$\sigma\ge0$ and $b\ge0$ are constants, $v\mapsto\gamma(v)$ is a
nondecreasing\vspace*{1pt} continuous function on $[0,1]$ such that $0\le
\gamma(v)\le1$ for all $0\le v\le1$ and $z^2\nu(dz)$ is a finite
measure on $(0,1]$. Let $\{B(ds,du)\}$ be a white noise on
$(0,\infty)\times(0,1]$ based on $ds\,du$, and let $\{M(ds,dz,du)\}$ be a
Poisson random measure on $(0,\infty)\times(0,1]^2$ with intensity
$ds\,\nu(dz)\,du$. We show that for any $v\in[0,1]$ there is a pathwise
unique solution $X(v) = \{X_t(v)\dvtx t\ge0\}$ to the equation
%
%
\begin{eqnarray}\label{1.6}
X_t(v) &=& v + \sigma\int_0^t\int_0^1 \bigl[1_{\{u\le X_{s-}(v)\}} -
X_{s-}(v)\bigr] B(ds,du) \nonumber\\
&&{}
+ b\int_0^t [\gamma(v)-X_{s-}(v)] \,ds \\
&&{}
+ \int_0^t\int_0^1\int_0^1z\bigl[1_{\{u\le X_{s-}(v)\}} - X_{s-}(v)\bigr]
M(ds,dz,du).\nonumber
\end{eqnarray}
Clearly, the above equation unifies and generalizes the flows
described\break
by~(\ref{1.1}), (\ref{1.2}) and (\ref{1.3}). Here it is essential to use
the white noise as the diffusion driving force. We show there is a
version of the random field $\{X_t(v)\dvtx\allowbreak t\ge0, 0\le v\le1\}$ with the
following properties:
\begin{longlist}
\item for each $v\in[0,1]$, $t\mapsto X_t(v)$ is c\`{a}dl\`{a}g on
$[0,\infty)$ and solves (\ref{1.6});

\item for each $t\ge0$, $v\mapsto X_t(v)$ is nondecreasing
and c\`{a}dl\`{a}g on $[0,1]$ with $X_t(0)\ge0$ and $X_t(1)\le1$.
\end{longlist}
We refer to $\{X_t(v)\dvtx t\ge0, 0\le v\le1\}$ as a \textit{generalized
Fleming--Viot flow} following Bertoin and Le Gall (\citeyear{BeL03},
\citeyear{BeL05}, \citeyear{BeL06}). In
particular, our result gives the strong existence of the flows associated
with the coalescents with multiple collisions. The study of this flow is
more involved than the one defined by (\ref{1.5}) as the path-valued
process $\{X(v)\dvtx0\le v\le1\}$ does not have independent increments.
However, we shall see it is still an inhomogeneous Markov process. {From}
the random field $\{X_t(v)\dvtx t\ge0, 0\le v\le1\}$ we can define a
c\`{a}dl\`{a}g sub-probability-valued process $\{X_t\dvtx t\ge0\}$ on $[0,1]$,
which is a counterpart of the generalized Fleming--Viot process of
\citet{BeL03}. We prove two scaling limit theorems for the
generalized Fleming--Viot processes, which lead to a special form of the
immigration superprocess defined from (\ref{1.5}). {From} the theorems we
derive easily a generalization of the limit theorem for the finite point
motions in \citet{BeL06}.

The techniques of this paper are mainly based on the strong solutions
of~(\ref{1.5}) and (\ref{1.6}), which are different from those of
Bertoin and Le Gall (\citeyear{BeL05}, \citeyear{BeL06}). In
Section \ref{sec2} we give some general
results for the pathwise uniqueness, comparison property and
existence of nonnegative strong solutions of stochastic equations
driven by white noises and Poisson random measures. Those extend the
results in \citet{FuL10} and provide the basis for the
investigation of the strong solution flows of (\ref{1.5}) and
(\ref{1.6}). They should also be of interest on their own right. In
Section \ref{sec3} we study the flows of CBI-processes and their associated
immigration superprocesses. The generalized Fleming--Viot flows are
discussed in Section \ref{sec4}. Finally, we prove the scaling limit theorems
in Section \ref{sec5}.

\subsection*{Notation}
For a measure $\mu$ and a function $f$ on a measurable space
$(E,\mcr{E})$ write $\langle\mu,f\rangle= \int_E f\,d\mu$ if the integral
exists. For any $a\ge0$ let $M[0,a]$ be the set of finite measures on
$[0,a]$ endowed with the topology of weak convergence. Let $M_1[0,a]$
be the subspace of $M[0,a]$ consisting of sub-probability measures. Let
$B[0,a]$ be the Banach space of bounded Borel functions on $[0,a]$
endowed with the supremum norm \mbox{$\|\cdot\|$}, and let $C[0,a]$ denote its
subspace of continuous functions. We use $B[0,a]^+$ and $C[0,a]^+$ to
denote the subclasses of nonnegative elements. Throughout this paper,
we make the conventions
\[
\int_a^b = \int_{(a,b]}
\quad\mbox{and}\quad
\int_a^\infty= \int_{(a,\infty)}
\]
for any $b\ge a\ge0$. Given a function $f$ defined on a subset of
$\mbb{R}$, we write
\[
\Delta_zf(x) = f(x+z) - f(x)
\quad\mbox{and}\quad
D_zf(x) = \Delta_zf(x) - f^\prime(x)z
\]
for $x,z\in\mbb{R}$ if the right-hand side is meaningful. Let
$\lambda$ denote the Lebesgue measure on $[0,\infty)$.


\section{Strong solutions of stochastic equations}\label{sec2}

In this section, we prove some results on stochastic equations of
one-dimensional processes driven by white noises and Poisson random
measures. The results extend those of \citet{FuL10}. Since our
aim is
to apply the results to the generalized Fleming--Viot flows and the flows
of CBI-processes, we only discuss equations of nonnegative processes.
However, the arguments can be modified to deal with general
one-dimensional equations.

Let $E$, $U_0$ and $U_1$ be
separable topological spaces whose topologies can be defined by
complete metrics. Suppose that $\pi(dz)$, $\mu_0(du)$ and
$\mu_1(du)$ are $\sigma$-fi\-nite Borel measures on $E$, $U_0$ and
$U_1$, respectively. We say the parameters~$(\sigma,b,g_0,g_1)$ are
\textit{admissible} if:
\begin{itemize}

\item$x\mapsto b(x)$ is a continuous function on $\mbb{R}_+$
satisfying $b(0)\ge0$;\vadjust{\goodbreak}

\item$(x,u)\mapsto\sigma(x,u)$ is a Borel function on
$\mbb{R}_+\times E$ satisfying $\sigma(0,u)=0$ for \mbox{$u\in
E$};

\item$(x,u) \mapsto g_0(x,u)$ is a Borel function on
$\mbb{R}_+\times U_0$ satisfying $g_0(0,u)=0$ and $g_0(x,u)+x\ge
0$ for $x>0$ and $u\in U_0$;

\item$(x,u) \mapsto g_1(x,u)$ is a Borel function on
$\mbb{R}_+\times U_1$ satisfying $g_1(x,u)+x\ge0$ for $x\ge0$
and $u\in U_1$.

\end{itemize}
Let $\{W(ds,du)\}$ be a white noise on $(0,\infty)\times E$ with
intensity $ds\,\pi(dz)$. Let $\{N_0(ds,du)\}$ and $\{N_1(ds,du)\}$ be
Poisson random measures on $(0,\infty)\times U_0$ and $(0,\infty
)\times
U_1$ with intensities $ds\,\mu_0(du)$ and $ds\,\mu_1(du)$, respectively.
Suppose that $\{W(ds,du)\}$, $\{N_0(ds,du)\}$ and $\{N_1(ds,du)\}$ are
defined on some complete probability space $(\Omega,\mcr{F},\mbf{P})$
and are independent of each other. Let $\{\tilde{N}_0(ds,du)\}$ denote
the compensated measure of $\{N_0(ds,du)\}$. A nonnegative c\`{a}dl\`{a}g
process $\{x(t)\dvtx t\ge0\}$ is called a \textit{solution} of
%
%
\begin{eqnarray}\label{2.1}
x(t)
&=&
x(0) + \int_0^t\int_{E} \sigma(x(s-),u)W(ds,du) \nonumber\\
&&{}
+ \int_0^t b(x(s-))\,ds + \int_0^t\int_{U_0} g_0(x(s-),u)\tilde{N}_0(ds,du)
\\
&&{}
+ \int_0^t\int_{U_1} g_1(x(s-),u) N_1(ds,du),\nonumber
\end{eqnarray}
if it satisfies the stochastic equation almost surely for every $t\ge0$.
We say $\{x(t)\dvtx t\ge0\}$ is a \textit{strong solution} if, in addition,
it is adapted to the augmented natural filtration generated by
$\{W(ds,du)\}$, $\{N_0(ds,du)\}$ and $\{N_1(ds,du)\}$ [see, e.g., Situ
(\citeyear{Sit05}), page 76]. Since $x(s-) \neq x(s)$ for at most
countably many $s\ge
0$, we can also use $x(s)$ instead of $x(s-)$ in the integrals with
respect to $W(ds,du)$ and $ds$ on the right-hand side of (\ref{2.1}). For
the convenience of the statements of the results, we write $b(x) =
b_1(x)-b_2(x)$, where $x\mapsto b_1(x)$ is continuous, and $x\mapsto
b_2(x)$ is continuous and nondecreasing. Let us formulate the following
conditions:
\begin{longlist}[(2.a)]
\item[(2.a)]
there is a constant $K\ge0$ so that
\[
b(x) + \int_{U_1} |g_1(x,u)| \mu_1(du) \le K(1+x)
\]
for every $x\ge0$;

\item[(2.b)] there is a nondecreasing function $x\mapsto L(x)$
on $\mbb{R}_+$ and a Borel function $(x,u)\mapsto
\bar{g}_0(x,u)$ on $\mbb{R}_+\times U_0$ so that $\sup_{0\le
y\le x}|g_0(y,u)|\le\bar{g}_0(x,u)$ and
\[
\int_E
\sigma(x,u)^2 \pi(du) + \int_{U_0} [\bar{g}_0(x,u)\land
\bar{g}_0(x,u)^2] \mu_0(du)\le L(x)
\]
for every $x\ge0$;

\item[(2.c)]
for each $m\ge1$ there is a nondecreasing concave function
$z\mapsto r_m(z)$ on $\mbb{R}_+$ such that $\int_{0+} r_m(z)^{-1}
\,dz=\infty$ and
\[
|b_1(x)-b_1(y)| + \int_{U_1} |g_1(x,u) - g_1(y,u)|\mu_1(du)
\le
r_m(|x-y|)
\]
for every $0\le x,y\le m$;

\item[(2.d)] for each $m\ge1$ there is a nonnegative
nondecreasing function $z\mapsto\rho_m(z)$ on $\mbb{R}_+$ so
that $\int_{0+} \rho_m(z)^{-2}\,dz = \infty$,
\[
\int_E
|\sigma(x,u)-\sigma(y,u)|^2 \pi(du) \le\rho_m(|x-y|)^2
\]
and
\[
\int_{U_0}\mu_0(du)\int_0^1 \frac{l_0(x,y,u)^2 (1-t)
1_{\{|l_0(x,y,u)|\le n\}}} {\rho_m(|(x-y) +
tl_0(x,y,u)|)^2} \,dt\le c(m,n)
\]
for every $n\ge1$
and $0\le x,y\le m$, where $l_0(x,y,u) = g_0(x,u) - g_0(y,u)$ and
$c(m,n)\ge0$ is a constant.
\end{longlist}
\begin{theorem}\label{t2.1} Suppose that $(\sigma,b,g_0,g_1)$ are admissible
parameters satisfying conditions \textup{(2.a)--(2.d)}. Then the pathwise
uniqueness of solutions holds for~(\ref{2.1}).
\end{theorem}
\begin{pf}
We first fix the integer $m\ge1$. Let $a_0=1$ and choose
$a_k\to0$ decreasingly\vspace*{1pt} so that $\int_{a_k} ^{a_{k-1}}
\rho_m(z)^{-2}\,dz = k$ for $k\ge1$. Let $x\mapsto\psi_k(x)$ be a~%
nonnegative continuous function on $\mbb{R}$ which has\vspace*{1pt} support in
$(a_k, a_{k-1})$ and satisfies $\int_{a_k}^{a_{k-1}} \psi_k(x) \,dx
=1$ and $0\le\psi_k(x) \le2 k^{-1}\rho_m(x)^{-2}$ for $a_k< x<
a_{k-1}$. For each $k\ge1$ we define the nonnegative and twice
continuously differentiable function
%
%
\begin{equation}\label{2.2}
\phi_k(z) = \int_0^{|z|}dy \int_0^y\psi_k(x)\,dx,\qquad z\in\mbb{R}.
\end{equation}
It is easy to see that $\phi_k(z)\to|z|$ nondecreasingly as $k\to
\infty$ and $0\le\phi_k^\prime(z) \le1$ for $z\ge0$ and $-1\le
\phi_k^\prime(z)\le0$ for $z\le0$. By condition (2.d) and the
choice of $x\mapsto\psi_k(x)$,
%
%
\begin{eqnarray}\label{2.3}
&&\phi_k^{\prime\prime}(x-y)\int_E |\sigma(x,u)-\sigma(y,u)|^2
\pi(du) \nonumber\\[-8pt]\\[-8pt]
&&\qquad
\le\psi_k(|x-y|)\rho_m(|x-y|)^2 \le\frac{2}{k}\nonumber
\end{eqnarray}
for $0\le x,y\le m$. Then the left-hand side tends to zero uniformly
in $0\le x,y\le m$ as $k\to\infty$. For $h,\zeta\in\mbb{R}$, by
Taylor's expansion we have
\[
D_h\phi_k(\zeta)
=
\int_0^1h^2\phi_k^{\prime\prime}(\zeta+th)(1-t)\,dt
=
\int_0^1h^2\psi_k(|\zeta+th|)(1-t)\,dt.
\]
It follows that
%
%
\begin{equation}\label{2.4}
D_h\phi_k(\zeta)
\le
\frac{2}{k}\int_0^1 h^2\rho_m(|\zeta+th|)^{-2}(1-t)\,dt.
\end{equation}
Observe also that
%
%
\begin{equation}\label{2.5}
D_h\phi_k(\zeta)
=
\Delta_h\phi_k(\zeta) - \phi_k^\prime(\zeta)h \le2|h|.
\end{equation}
For $0\le x,y\le m$ and $n\ge1$ we can use (\ref{2.4}) and (\ref
{2.5}) to
get
%
%
\begin{eqnarray}\label{2.6}
&&\int_{U_0} D_{l_0(x,y,u)}\phi_k(x-y) \mu_0(du) \nonumber\\
&&\qquad
\le\frac{2}{k}\int_{U_0}\mu_0(du)\int_0^1 \frac{l_0(x,y,u)^2(1-t)
1_{\{|l_0(x,y,u)|\le n\}}}{\rho_m(|(x-y) +
tl_0(x,y,u)|)^2} \,dt \nonumber\\[-8pt]\\[-8pt]
&&\qquad\quad{}
+ 2\int_{U_0} |l_0(x,y,u)|1_{\{|l_0(x,y,u)|>n\}} \mu_0(du) \nonumber\\
&&\qquad
\le\frac{2}{k}c(m,n) + 4\int_{U_0} \bar{g}_0(m,u)
1_{\{\bar{g}_0(m,u)>n/2\}} \mu_0(du).\nonumber
\end{eqnarray}
By conditions (2.b), (2.d) one sees the right-hand side tends to zero
uniformly in $0\le x,y\le m$ as $k\to\infty$. Then the pathwise
uniqueness for (\ref{2.1}) follows by a trivial modification of
Theorem 3.1 in \citet{FuL10}.
\end{pf}

The key difference between the above theorem and Theorems 3.2
and 3.3 of \citet{FuL10} is that here we do not assume $x\mapsto
g_0(x,u)$ is nondecreasing. This is essential for the applications
to stochastic equations like~(\ref{1.6}).

\begin{theorem}\label{t2.2} Let $(\sigma,b^\prime,g_0,g_1^\prime)$ and
$(\sigma,b^{\prime\prime},g_0,g_1^{\prime\prime})$ be two sets of
admissible parameters satisfying conditions \textup{(2.a)--(2.d)}. In
addition, assume that:
\begin{longlist}
\item for every $u\in U_1$, $x \mapsto x + g_1^\prime(x,u)$
or $x \mapsto x + g_1^{\prime\prime}(x,u)$ is nondecreasing;

\item$b^\prime(x)\le b^{\prime\prime}(x)$ and
$g_1^\prime(x,u)\le g_1^{\prime\prime}(x,u)$ for every $x\ge0$
and $u\in U_1$.
\end{longlist}
Suppose that $\{x^\prime(t)\dvtx t\ge0\}$ is a solution of (\ref{2.1})
with $(b,g_1) = (b^\prime,g_1^\prime)$, and $\{x^{\prime\prime}(t)\dvtx
t\ge0\}$ is a solution of the equation with $(b,g_1) =
(b^{\prime\prime},g_1^{\prime\prime})$. If $x^\prime(0)\le
x^{\prime\prime}(0)$, then $\mbf{P}\{x^\prime(t)\le
x^{\prime\prime}(t)$ for all $t\ge0\} = 1$.
\end{theorem}
\begin{pf}
Let $\zeta(t) = x^\prime(t) - x^{\prime\prime}(t)$ for
$t\ge
0$. Let $x\mapsto\psi_k(x)$ be defined as in the proof of
Theorem \ref{t2.1}. Instead of (\ref{2.2}), for each $k\ge1$ we now
define
%
%
\begin{equation}\label{2.7}
\phi_k(z) = \int_0^z dy \int_0^y\psi_k(x)\,dx,\qquad z\in\mbb{R}.
\end{equation}
Then $\phi_k(z) \to z^+ := 0\vee z$ nondecreasingly as $k\to
\infty$. Let
\[
l_0(t,u) = g_0(x^\prime(t),u) - g_0(x^{\prime\prime}(t),u),\qquad
t\ge0, u\in U_0,\vadjust{\goodbreak}
\]
and
\[
l_1(t,u) = g_1^\prime(x^\prime(t),u) - g_1^{\prime\prime}
(x^{\prime\prime}(t),u),\qquad t\ge0, u\in U_1.
\]
For $\zeta(s-)\le0$ we have $\phi_k(\zeta(s-)) =
\phi_k^\prime(\zeta(s-)) = 0$. Since $x \mapsto x + f(x,u)$ is
nondecreasing for $f = g_1^\prime$ or $g_1^{\prime\prime}$, for
$\zeta(s-) = x^\prime(s-) - x^{\prime\prime}(s-)\le0$ we also have
\begin{eqnarray*}
\zeta(s-)+l_1(s-,u)
&=&
x^\prime(s-) - x^{\prime\prime}(s-) + g_1^\prime(x^\prime(s-),u) -
g_1^{\prime\prime} (x^{\prime\prime}(s-),u) \\
&\le&
x^\prime(s-) - x^{\prime\prime}(s-) + f(x^\prime(s-),u) -
f(x^{\prime\prime}(s-),u)
\le0.
\end{eqnarray*}
The latter implies
\[
\Delta_{l_1(s-,u)} \phi_k(\zeta(s-))
=
\phi_k\bigl(\zeta(s-)+l_1(s-,u)\bigr) - \phi_k(\zeta(s-)) = 0.
\]
By It\^{o}'s formula we have
%
%
\begin{eqnarray}\label{2.8}
\phi_k(\zeta(t))
&=&
\phi_k(\zeta(0)) + \frac{1}{2}\int_0^t\,ds \int_E
\phi_k^{\prime\prime}(\zeta(s-)) [\sigma(x^\prime(s-),u) \nonumber\\
&&\hspace*{153pt}{}
- \sigma(x^{\prime\prime}(s-),u)]^2 \pi(du) \nonumber\\
&&{}
+ \int_0^t \phi_k^\prime(\zeta(s-)) [b^\prime(x^\prime(s-))
\nonumber\\[-8pt]\\[-8pt]
&&\hspace*{76pt}{}
- b^{\prime\prime}(x^{\prime\prime}(s-))] 1_{\{\zeta(s-)>0\}}\,ds
\nonumber\\
&&{}
+ \int_0^t ds \int_{U_1} \Delta_{l_1(s-,u)} \phi_k(\zeta(s-))
1_{\{\zeta(s-)>0\}}\mu_1(du) \nonumber\\
&&{}
+ \int_0^t ds \int_{U_0} D_{l_0(s-,u)} \phi_k(\zeta(s-)) \mu_0(du) +
M_m(t),\nonumber
\end{eqnarray}
where
\begin{eqnarray*}
M_m(t)
&=&
\int_0^t\int_E\phi_k^\prime(\zeta(s-)) [\sigma(x^\prime(s-),u)\\
&&\hspace*{77pt}{}
- \sigma(x^{\prime\prime}(s-),u)]W(ds,du) \\
&&{}
+ \int_0^t\int_{U_1} \Delta_{l_1(s-,u)} \phi_k(\zeta(s-))
\tilde{N}_1(ds,du) \\
&&{}
+ \int_0^t\int_{U_0} \Delta_{l_0(s-,u)}\phi_k(\zeta(s-))
\tilde{N}_0(ds,du).
\end{eqnarray*}
Let $\tau_m = \inf\{t\ge0\dvtx x^\prime(t)\ge m$ or $x^{\prime\prime
}(t)\ge
m\}$ for $m\ge1$. Under conditions~(2.b), (2.c) it is easy to show that
$\{M_m(t\land\tau_m)\}$ is a martingale. Recall that $b^\prime(x)\le
b^{\prime\prime}(x)$ and $b^\prime(x) = b^\prime_1(x) - b^\prime_2(x)$
for a nondecreasing\vadjust{\goodbreak} function $x\mapsto b^\prime_2(x)$.
Then under the restriction $\zeta(s-)>0$ we have
\begin{eqnarray*}
&&\phi_k^\prime(\zeta(s-))[b^\prime(x^\prime(s-)) -
b^{\prime\prime}(x^{\prime\prime}(s-))] \\
&&\qquad
\le\phi_k^\prime(\zeta(s-))[b^\prime(x^\prime(s-)) -
b^\prime(x^{\prime\prime}(s-))] \\
&&\qquad
\le\phi_k^\prime(\zeta(s-))[b_1^\prime(x^\prime(s-)) -
b_1^\prime(x^{\prime\prime}(s-))] \\
&&\qquad
\le|b_1^\prime(x^\prime(s-))-b_1^\prime(x^{\prime\prime}(s-))|
\end{eqnarray*}
and
\begin{eqnarray*}
&&\Delta_{l_1(s-,u)} \phi_k(\zeta(s-)) \\
&&\qquad
= \phi_k\bigl(\zeta(s-)+g_1^\prime(x^\prime(s-),u) -
g_1^{\prime\prime}(x^{\prime\prime}(s-),u)\bigr) - \phi_k(\zeta(s-))
\\
&&\qquad
\le\phi_k\bigl(\zeta(s-)+g_1^\prime(x^\prime(s-),u) -
g_1^\prime(x^{\prime\prime}(s-),u)\bigr) - \phi_k(\zeta(s-)) \\
&&\qquad
\le|g_1^\prime(x^\prime(s-),u) - g_1^\prime(x^{\prime\prime}(s-),u)|.
\end{eqnarray*}
The estimates (\ref{2.3}) and (\ref{2.6}) are still valid. If
$x^\prime(0)\le x^{\prime\prime}(0)$, we can take the expectation in
(\ref{2.8}) and let $k\to\infty$ to get
\begin{eqnarray*}
\mbf{E}[\zeta(t\land\tau_m)^+]
&\le&
\mbf{E}\biggl[\int_0^{t\land\tau_m} r_m(|\zeta(s-)|)
1_{\{\zeta(s-)>0\}} \,ds\biggr] \\
&\le&
\int_0^t r_m\bigl(\mbf{E}[\zeta(s\land\tau_m)^+]\bigr)\,ds,
\end{eqnarray*}
where the second inequality holds by the concaveness of $z\mapsto
r_m(z)$. Then $\mbf{E}[\zeta(t\land\tau_m)^+] =0$ for all $t\ge0$.
Since $\tau_m\to\infty$ as $m\to\infty$, we get the desired comparison
property.
\end{pf}

We say the \textit{comparison property} of solutions holds for
(\ref{2.1}) if for any two solutions $\{x_1(t)\dvtx t\ge0\}$ and
$\{x_2(t)\dvtx
t\ge0\}$ satisfying $x_1(0)\le x_2(0)$ we have $\mbf{P} \{x_1(t)\le
x_2(t)$ for all $t\ge0\} = 1$. {From} Theorem \ref{t2.2} we get the
following:
\begin{theorem}\label{t2.3} Let $(\sigma,b,g_0,g_1)$ be admissible parameters
satisfying conditions \textup{(2.a)--(2.d)}. In addition,
assume that for
every $u\in U_1$ the function $x \mapsto x + g_1(x,u)$ is nondecreasing.
Then the comparison property holds for the solutions of (\ref{2.1}).
\end{theorem}

The monotonicity assumption on the function $x \mapsto x + g_1(x,u)$ in
Theorem \ref{t2.3} is natural. To see this, suppose that $\{x_1(t)\}$ and
$\{x_2(t)\}$ are two solutions of (\ref{2.1}) and $\{(s_i,u_i)\dvtx i\ge
1\}$
is the set of atoms of $\{N_1(ds,du)\}$. The assumption guarantees that
$x_1(s_i-)\le x_2(s_i-)$ implies
\begin{eqnarray*}
x_1(s_i)
&=&
x_1(s_i-) + g_1(x_1(s_i-),u_i) \\
&\le&
x_2(s_i-) + g_1(x_2(s_i-),u_i) = x_2(s_i).
\end{eqnarray*}
A similar explanation can be given to Theorem \ref{t2.2}. In some
applications the kernel $x \mapsto g_0(x,u)$ may be nondecreasing. When
this is true, we can replace~(2.d) by the following simpler
condition:
\begin{longlist}[(2.e)]
\item[(2.e)]
For each $u\in U_0$ the function $x \mapsto g_0(x,u)$ is
nondecreasing, and for each $m\ge1$ there is a nonnegative and
nondecreasing function $z\mapsto\rho_m(z)$ on~$\mbb{R}_+$ so that
$\int_{0+} \rho_m(z)^{-2}\,dz = \infty$ and
\begin{eqnarray*}
&&\int_E |\sigma(x,u)-\sigma(y,u)|^2 \pi(du)
+ \int_{U_0}
|l_0(x,y,u)|\land|l_0(x,y,u)|^2 \mu_0(du) \\
&&\qquad
\le\rho_m(|x-y|)^2
\end{eqnarray*}
for all $0\le x,y\le m$, where $l_0(x,y,u) = g_0(x,u) -
g_0(y,u)$.
\end{longlist}
\begin{prop}\label{p2.4} Let $(\sigma,b,g_0,g_1)$ be admissible
parameters. If \textup{(2.e)} holds, then \textup{(2.d)} holds.
\end{prop}
\begin{pf}
Since $x \mapsto g_0(x,u)$ is nondecreasing, it is not hard to
see $|(x-y) + tl_0(x,y,u)|\ge|x-y|$. By condition \textup{(2.e)} and the
monotonicity of $z\mapsto\rho(z)$ we have
\begin{eqnarray*}
&&\int_0^1dt\int_{U_0} \frac{(1-t)l_0(x,y,u)^2 1_{\{
|l_0(x,y,u)|\le
n\}}} {\rho_m(|(x-y) + tl_0(x,y,u)|)^2} \mu_0(du)
\\
&&\qquad
\le n\int_0^1dt\int_{U_0} \frac{[|l_0(x,y,u)|\land l_0(x,y,u)^2]}
{\rho_m(|x-y|)^2} \mu_0(du) \le n.
\end{eqnarray*}
Then condition \textup{(2.d)} is satisfied.
\end{pf}
\begin{theorem}\label{t2.5} Suppose that $(\sigma,b,g_0,g_1)$ are admissible
parameters satisfying conditions \textup{(2.a), (2.c), (2.e)}. Then
there is a unique
strong solution to (\ref{2.1}).
\end{theorem}
\begin{pf}
We first note that \textup{(2.b)} follows from \textup{(2.e)}. By
Proposition \ref{p2.4}, we also have \textup{(2.d)} from \textup
{(2.e)}. Let
$\{V_n\}$ be a nondecreasing sequence of Borel subsets of $U_0$ so that
$\bigcup_{n=1}^\infty V_n = U_0$ and $\mu_0(V_n)< \infty$ for every
$n\ge
1$. For $m,n\ge1$ one can use \textup{(2.e)} to see
\[
x\mapsto\beta_m(x) := \int_{U_0} [g_0(x,u) - g_0(x,u)\land m] \mu_0(du)
\]
and
\[
x\mapsto\gamma_{m,n}(x) := \int_{V_n} [g_0(x,u)\land m] \mu_0(du)
\]
are continuous nondecreasing functions. By the results for
continuous-type stochastic
equations as in Ikeda and Watanabe\vadjust{\goodbreak} [(\citeyear{IkW89}), page 169] one can
show there is a nonnegative weak solution to
\begin{eqnarray*}
x(t) &=& x(0) + \int_0^t\int_{E} \sigma\bigl(x(s)\land m,u\bigr) W(ds,du)
\\[-2pt]
&&{}
+ \int_0^t b_m\bigl(x(s)\land m\bigr) \,ds - \int_0^t \gamma_{m,n}\bigl(x(s)\land m\bigr) \,ds,
\end{eqnarray*}
where $b_m(x) = b(x) - \beta_m(x)$. The pathwise uniqueness holds for
the above equation by
Theorem \ref{t2.1}. Then it has a unique strong solution. Let $\{W_n\}$
be a nondecreasing sequence of Borel subsets of $U_1$ so that
$\bigcup_{n=1}^\infty W_n = U_1$ and $\mu_1(W_n)< \infty$ for every
$n\ge
1$. Following the proof of Proposition 2.2 of \citet{FuL10} one can
show there is a unique strong solution $\{x_{m,n}(t)\dvtx t\ge0\}$
to\looseness=-1
\begin{eqnarray*}
x(t)
&=&
x(0) + \int_0^t\int_{E} \sigma\bigl(x(s-)\land m,u\bigr)W(ds,du) \\[-2pt]
&&{}
+ \int_0^t b_m\bigl(x(s-)\land m\bigr)\,ds - \int_0^t \gamma_{m,n}\bigl(x(s)\land m\bigr) \,ds
\\[-2pt]
&&{}
+ \int_0^t\int_{V_n} \bigl[g_0\bigl(x(s-)\land m,u\bigr)\land m\bigr] N_0(ds,du) \\[-2pt]
&&{}
+ \int_0^t\int_{W_n} \bigl[g_1\bigl(x(s-)\land m,u\bigr)\land m\bigr] N_1(ds,du).
\end{eqnarray*}\looseness=0
We can rewrite the above equation into
\begin{eqnarray*}
x(t)
&=&
x(0) + \int_0^t\int_{E} \sigma\bigl(x(s-)\land m,u\bigr)W(ds,du) \\[-2pt]
&&{}
+ \int_0^t b_m\bigl(x(s-)\land m\bigr)\,ds \\[-2pt]
&&{}
+ \int_0^t\int_{V_n} \bigl[g_0\bigl(x(s-)\land m,u\bigr)\land m\bigr] \tilde{N}_0(ds,du)
\\[-2pt]
&&{}
+ \int_0^t\int_{W_n} \bigl[g_1\bigl(x(s-)\land m,u\bigr)\land m\bigr] N_1(ds,du).
\end{eqnarray*}
As in the proof of Lemma 4.3 of \citet{FuL10} one can see the sequence
$\{x_{m,n}(t)\dvtx t\ge0\}$, $n=1,2, \ldots,$ is tight in
$D([0,\infty),\mbb{R}_+)$. Following the proof of Theorem 4.4 of
\citet{FuL10} it is easy to show that any weak limit point $\{
x_m(t)\dvtx t\ge
0\}$ of the sequence is a nonnegative weak solution to
%
%
\begin{eqnarray}\label{2.9}
x(t)
&=&
x(0) + \int_0^t\int_{E} \sigma\bigl(x(s-)\land m,u\bigr)W(ds,du) \nonumber\\[-2pt]
&&{}
+ \int_0^t b_m\bigl(x(s-)\land m\bigr)\,ds \nonumber\\[-9pt]\\[-9pt]
&&{}
+ \int_0^t\int_{U_0} \bigl[g_0\bigl(x(s-)\land m,u\bigr)\land m\bigr] \tilde{N}_0(ds,du)
\nonumber\\[-2pt]
&&{}
+ \int_0^t\int_{U_1} \bigl[g_1\bigl(x(s-)\land m,u\bigr)\land m\bigr] N_1(ds,du).\nonumber
\end{eqnarray}
By Theorem \ref{t2.1} the pathwise uniqueness holds for (\ref{2.9}), so
the equation has a unique strong solution [see, e.g., Situ (\citeyear
{Sit05}), page 104].
Then the result follows by a simple modification of the proof of
Proposition 2.4 of \citet{FuL10}. See \citet{ElM90}
and Kurtz (\citeyear{Kur07}, \citeyear{Kur10}) for the general theory
of stochastic
equations driven by white noises and Poisson random measures.
\end{pf}


\section{Stochastic flows of CBI-processes}\label{sec3}

In this section, we give the constructions and characterizations of
the flow of CBI-processes and the associated immigration
superprocess. Suppose that $\sigma\ge0$ and $b$ are constants, and
$(u\land u^2)m(du)$ is a finite measure on $(0,\infty)$. Let $\phi$
be a function given by
%
%
\begin{equation}\label{3.1}
\phi(z) = bz + \frac{1}{2}\sigma^2z^2 + \int_0^\infty
(e^{-zu}-1+zu)m(du),\qquad z\ge0.
\end{equation}
A Markov process with state space $\mbb{R}_+ := [0,\infty)$ is called a
\textit{CB-process} with branching mechanism $\phi$ if it has transition
semigroup $(p_t)_{t\ge0}$ given by
%
%
\begin{equation}\label{3.2}
\int_{\mbb{R}_+} e^{-\lambda y} p_t(x,dy) = e^{-xv_t(\lambda)},\qquad
\lambda\ge0,
\end{equation}
where $(t,\lambda)\mapsto v_t(\lambda)$ is the unique locally
bounded nonnegative solution of
\[
\frac{d}{dt}v_t(\lambda) = - \phi(v_t(\lambda)),\qquad v_0(\lambda)
= \lambda,\qquad t\ge0.
\]
Given any $\beta\ge0$ we can also define a transition semigroup
$(q_t)_{t\ge0}$ on $\mbb{R}_+$ by
%
%
\begin{equation}\label{3.3}
\int_{\mbb{R}_+} e^{-\lambda y} q_t(x,dy)
=
\exp\biggl\{-xv_t(\lambda)-\int_0^t\beta v_s(\lambda) \,ds\biggr\}.
\end{equation}
A nonnegative real-valued Markov process with transition semigroup
$(q_t)_{t\ge0}$ is called a \textit{CBI-process} with branching mechanism
$\phi$ and immigration rate~$\beta$. It is easy to see that both
$(p_t)_{t\ge0}$ and $(q_t)_{t\ge0}$ are Feller semigroups. See, for example,
\citet{KaW71} and Li (\citeyear{Li10}), Chapter 3.\vspace*{1pt}

Let $\{W(ds,du)\}$ be a white noise on $(0,\infty)^2$ based on the
Lebesgue measure $ds\,du$, and let $\{N(ds,dz,du)\}$ be Poisson random
measure on $(0,\infty)^3$ with intensity $ds\,m(dz)\,du$. Let
$\{\tilde{N}(ds,dz,du)\}$ be the compensated measure of
$\{N(ds,dz,du)\}$.
\begin{theorem}\label{t3.1} There is a unique nonnegative strong solution
of the stochastic equation
\begin{eqnarray*}
Y_t &=& Y_0 + \sigma\int_0^t\int_0^{Y_{s-}} W(ds,du) + \int_0^t
(\beta-bY_{s-}) \,ds \\
&&{}
+ \int_0^t \int_0^\infty\int_0^{Y_{s-}} z \tilde{N}(ds,dz,du).
\end{eqnarray*}
Moreover, the solution $\{Y_t\dvtx t\ge0\}$ is a CBI-process with
branching mechanism $\phi$ and immigration rate $\beta$.
\end{theorem}
\begin{pf}
The existence and uniqueness of the strong solution follows
by an
application of Theorem \ref{t2.5} [see also \citet{DaL06}]. Using
It\^{o}'s formula one can see that $\{Y_t(v)\dvtx t\ge0\}$ solves the
martingale problem associated with the generator $L$ defined by
%
%
\begin{equation}\label{3.4}
Lf(x) = \frac{1}{2}\sigma^2xf^{\prime\prime}(x) + (\beta- bx)
f^\prime(x) + x \int_0^\infty D_zf(x) m(dz).
\end{equation}
Then it is a CBI-process with branching mechanism $\phi$ and immigration
rate~$\beta$ [see \citet{KaW71} and Li (\citeyear{Li10}), Section 9.5].
\end{pf}

Let $v\mapsto\gamma(v)$ be a nonnegative and nondecreasing continuous
function on $[0,\infty)$. We denote by $\gamma(dv)$ the Radon measure on
$[0,\infty)$ so that $\gamma([0,v]) = \gamma(v)$ for $v\ge0$. By
Theorem \ref{t3.1} for each $v\ge0$ there is a pathwise unique
nonnegative solution $Y(v) = \{Y_t(v)\dvtx t\ge0\}$ to the stochastic
equation
%
%
\begin{eqnarray}\label{3.5}
Y_t(v) &=& v + \sigma\int_0^t\int_0^{Y_{s-}(v)} W(ds,du) +
\int_0^t [\gamma(v)-bY_{s-}(v)] \,ds \nonumber\\[-8pt]\\[-8pt]
&&{}
+ \int_0^t \int_0^\infty\int_0^{Y_{s-}(v)} z
\tilde{N}(ds,dz,du).\nonumber
\end{eqnarray}
\begin{theorem}\label{t3.2} For any $v_2\ge v_1\ge0$ we have
$\mbf{P}\{Y_t(v_2)\ge Y_t(v_1)$ for all $t\ge0\} = 1$ and $\{Y_t(v_2) -
Y_t(v_1)\dvtx t\ge0\}$ is a CBI-process with branching mechanism $\phi$ and
immigration rate $\beta:= \gamma(v_2)-\gamma(v_1)\ge0$.
\end{theorem}
\begin{pf}
The comparison property follows by applying
Theorem \ref{t2.2} and Proposition \ref{p2.4} to (\ref{3.5}). Let
$Z_t = Y_t(v_2) - Y_t(v_1)$ for $t\ge0$. {From} (\ref{3.5}) we
have\looseness=-1
%
%
\begin{eqnarray}\label{3.6}
Z_t &=& v_2-v_1 + \sigma\int_0^t\int_{Y_{s-}(v_1)}^{Y_{s-}(v_2)}
W(ds,du) + \int_0^t (\beta-bZ_{s-}) \,ds \nonumber\\
&&{}
+ \int_0^t \int_0^\infty\int_{Y_{s-}(v_1)}^{Y_{s-}(v_2)}
z\tilde{N}(ds,dz,du) \nonumber\\[-8pt]\\[-8pt]
&=&
v_2-v_1 + \sigma\int_0^t\int_0^{Z_{s-}} W_1(ds,du) + \int_0^t
(\beta-bZ_{s-})\,ds \nonumber\\
&&{}
+ \int_0^t \int_0^\infty\int_0^{Z_{s-}}z
\tilde{N}_1(ds,dz,du),\nonumber
\end{eqnarray}\looseness=0
where
\[
W_1(ds,du) = W\bigl(ds,Y_{s-}(v_1)+du\bigr)
\]
is a white noise with intensity $ds\,du$, and
\[
N_1(ds,dz,du) = N\bigl(ds,dz,Y_{s-}(v_1)+du\bigr)
\]
is a Poisson random measure with intensity $ds\,m(dz)\,du$. That shows
$\{Z_t\dvtx t\ge0\}$ is a weak solution of (\ref{3.5}). Then it a
CBI-process with branching mechanism $\phi$ and immigration rate
$\beta$.
\end{pf}
\begin{theorem}\label{t3.3} Let $v_2\ge v_1\ge u_2\ge u_1\ge0$. Then
$\{Y_t(u_2) - Y_t(u_1)\dvtx t\ge0\}$ and $\{Y_t(v_2) - Y_t(v_1)\dvtx
t\ge
0\}$ are independent CBI-processes with immigration rates $\alpha:=
\gamma(u_2)-\gamma(u_1)$ and $\beta:= \gamma(v_2)-\gamma(v_1)$,
respectively.
\end{theorem}
\begin{pf}
Let $L_\alpha$ and $L_\beta$ denote the generators of the
CBI-processes with immigration rates $\alpha$ and $\beta$,
respectively. Let $X_t = Y_t(u_2) - Y_t(u_1)$ and $Z_t = Y_t(v_2) -
Y_t(v_1)$. For any $G\in C^2(\mbb{R}_+^2)$ one can use It\^{o}'s
formula to show
%
%
\begin{eqnarray}\label{3.7}
G(X_t,Z_t) &=& G(X_0,Z_0) + \int_0^t L_\alpha G(X_s,Z_s) \,ds
\nonumber\\[-8pt]\\[-8pt]
&&{}
+ \int_0^t L_\beta G(X_s,Z_s) \,ds + \mbox{local mart.},\nonumber
\end{eqnarray}
where $L_\alpha$ and $L_\beta$ act on the first and second
coordinates of $G$, respectively. Then $\{X_t\dvtx t\ge0\}$ and $\{
Z_t\dvtx
t\ge0\}$ are independent CBI-processes with immigration rates
$\alpha$ and $\beta$, respectively.
\end{pf}
\begin{prop}\label{p3.4} There is a locally bounded nonnegative
function $t\mapsto C(t)$ on $[0,\infty)$ so that
%
%
\begin{eqnarray}\label{3.8}
\mbf{E}\Bigl\{\sup_{0\le s\le t}[Y_s(v_2) - Y_s(v_1)]\Bigr\}
&\le&
C(t)\bigl\{(v_2-v_1) + [\gamma(v_2)-\gamma(v_1)] \nonumber\\[-8pt]\\[-8pt]
&&\hspace*{23.5pt}{}
+ \sqrt{v_2-v_1} + \sqrt{\gamma(v_2)-\gamma(v_1)}\bigr\}\nonumber
\end{eqnarray}
for $t\ge0$ and $v_2\ge v_1\ge0$.
\end{prop}
\begin{pf}
Let $Z_t = Y_t(v_2)-Y_t(v_1)$ for $t\ge0$. Taking the
expectation in~(\ref{3.6}) we have
\[
\mbf{E}(Z_t)
=
(v_2-v_1) + t[\gamma(v_2)-\gamma(v_1)] - b\int_0^t \mbf{E}(Z_s)
\,ds.\vadjust{\goodbreak}
\]
Solving the above integral equation gives
%
%
\begin{equation}\label{3.9}
\mbf{E}(Z_t) = (v_2-v_1)e^{-bt} + [\gamma(v_2)-\gamma(v_1)]
b^{-1}(1-e^{-bt})
\end{equation}
with $b^{-1}(1-e^{-bt}) = t$ for $b=0$ by convention. By (\ref{3.6}) and
Doob's martingale inequality,
\begin{eqnarray*}
\mbf{E}\Bigl\{\sup_{0\le s\le t}Z_s\Bigr\}
&\le&
(v_2-v_1) + 2\sigma\mbf{E}^{1/2}\biggl\{\biggl(\int_0^t
\int_{Y_{s-}(v_1)}^{Y_{s-}(v_2)} W(ds,du)\biggr)^2\biggr\} \\
&&{}
+ \int_0^t \{[\gamma(v_2)-\gamma(v_1)] + |b|\mbf{E}(Z_s)\} \,ds \\
&&{}
+ 2\mbf{E}^{1/2}\biggl\{\biggl(\int_0^t\int_0^1
\int_{Y_{s-}(v_1)}^{Y_{s-}(v_2)} z \tilde{N}(ds,dz,du)\biggr)^2
\biggr\} \\
&&{}
+ \mbf{E}\biggl[\int_0^t\int_1^\infty\int
_{Y_{s-}(v_1)}^{Y_{s-}(v_2)} z
N(ds,dz,du)\biggr] \\
&\le&
(v_2-v_1) + t[\gamma(v_2)-\gamma(v_1)] + 2\sigma\biggl[\int_0^t
\mbf{E}(Z_s) \,ds\biggr]^{1/2} \\
&&{}
+ 2\biggl[\int_0^1 z^2\nu(dz)\biggr]^{1/2} \biggl[\int_0^t
\mbf{E}(Z_s) \,ds\biggr]^{1/2} \\
&&{}
+ \biggl[|b|+\int_1^\infty z \nu(dz)\biggr] \int_0^t \mbf{E}(Z_s)\,ds.
\end{eqnarray*}
Then (\ref{3.8}) follows by (\ref{3.9}).
\end{pf}

Suppose that $(E,\rho)$ is a complete metric space. Let $F$ be a subset
of $[0,\infty)$ such that $0\in F$ and let $t\mapsto x(t)$ be a path from
$F$ to $E$. For any $\varepsilon>0$ the number of $\varepsilon
$-oscillations of
this path on $F$ is defined as
\begin{eqnarray*}
\mu(\varepsilon) &:=&\sup\{n\ge0\dvtx \mbox{there are $0=t_0< t_1<
\cdots< t_n\in F$} \\
&&\hspace*{19pt}
\mbox{so that $\rho(x(t_{i-1}),x(t_i))\ge\varepsilon$ for all $1\le
i\le n\}$}.
\end{eqnarray*}
If $F$ is dense in $[0,\infty)$, it is simple to show the limits $y(t) :=
\lim_{F\ni s\to t+} x(s)$ exist for all $t\ge0$ and constitute a
c\`{a}dl\`{a}g path $t\mapsto y(t)$ on $[0,\infty)$ if and only if
$t\mapsto
x(t)$ has at most a finite number of $\varepsilon$-oscillations on
$F\cap
[0,T]$ for every $\varepsilon>0$ and $T\ge0$.
\begin{lemma}\label{l3.5} Suppose that $(\Omega, \mcr{G}, \mcr{G}_t,
\mbf{P})$ is a filtered probability space and $\{X_t\dvtx t\ge0\}$ is a
$(\mcr{G}_t)$-Markov process with state space $(E,\mcr{E})$ and
transition semigroup $(P_{s,t})_{t\ge s}$. Suppose that $\rho$ is a
complete metric on $E$ so that:
\begin{longlist}
\item for $\varepsilon>0$ and $0\le s,t\le u$ we have
$\{\omega\in\Omega\dvtx \rho(X_s(\omega),X_t(\omega))<
\varepsilon\}\in\mcr{G}_u$;

\item for $\varepsilon>0$ and $x\in E$ we have $U_\varepsilon(x)
:= \{y\in E\dvtx \rho(x,y)< \varepsilon\}\in\mcr{E}$ and
%
%
\begin{equation}\label{3.10}
\alpha_\varepsilon(h) := \sup_{0\le t-s\le h}\sup_{x\in E}
P_{s,t}(x,U_\varepsilon(x)^c)\to0\qquad (h\to0).
\end{equation}
\end{longlist}
Then $\{X_t\dvtx t\ge0\}$ has a $\rho$-c\`{a}dl\`{a}g modification.
\end{lemma}
\begin{pf}
Let $F = \{0,r_1,r_2,\ldots\}$ be a countable dense subset of
$[0,\infty)$ and let $F_n = \{0,r_1,\ldots,r_n\}$. For $\varepsilon>0$ and
$a>0$ let $\nu^a(\varepsilon)$ and $\nu_n^a(\varepsilon)$ denote,
respectively, the numbers of $\varepsilon$-oscillations of $t\mapsto
X_t$ on $F\cap[0,a]$
and $F_n\cap[0,a]$. Then $\nu_n^a(\varepsilon)\to\nu^a(\varepsilon)$
increasingly as $n\to\infty$. Let $\tau_n^\varepsilon(0) = 0$ and for
$k\ge
0$ define
\[
\tau_n^\varepsilon(k+1) = \min\bigl\{t\in F_n\cap(\tau_n^\varepsilon(k),
\infty)\dvtx
\rho\bigl(X_{\tau_n^\varepsilon(k)},X_t\bigr)\ge\varepsilon\bigr\},
\]
if $\tau_n^\varepsilon(k)< \infty$ and $\tau_n^\varepsilon(k+1) = \infty
$ if
$\tau_n^\varepsilon(k) = \infty$. Since $F_n$ is discrete, for any
$a\ge0$
we have
\[
\{\tau_n^\varepsilon(k+1)\le a\}
=
\bigcup_{s<t\in F_n\cap[0,a]} \bigl(\{\tau_n^\varepsilon(k) = s\}\cap
\{\rho(X_s,X_t)\ge\varepsilon\}\bigr).
\]
Using property (i) and the above relation it is easy to see successively
that each $\tau_n^\varepsilon(k)$ is a stopping time. As in the proof of
Lemma 9.1 of Wentzell [(\citeyear{Wen81}), page~168] one can prove
$\mbf{P}\{\tau_n^\varepsilon(1)\le h\}\le2\alpha_{\varepsilon/2}(h)$ for
$\varepsilon>0$ and $h>0$. Since the strong Markov property of $\{
X_t\dvtx
t\ge
0\}$ holds at the discrete stopping times $\tau_n^\varepsilon(k)$,
$k=1,2,\ldots,$ one can inductively show
\[
\mbf{P}\{\nu_n^h(2\varepsilon)\ge k\}
\le
\mbf{P}\{\tau_n^\varepsilon(k)\le h\}
\le
[2\alpha_{\varepsilon/2}(h)]^k.
\]
It follows that
\[
\mbf{P}\{\nu^h(2\varepsilon)\ge k\}
=
\lim_{n\to\infty}\mbf{P}\{\nu_n^h(2\varepsilon)\ge k\}
\le
[2\alpha_{\varepsilon/2}(h)]^k.
\]
Choosing sufficiently small $h=h(\varepsilon)\in F\cap(0,\infty)$ so that
$\alpha_{\varepsilon/2}(h)< 1/2$ and letting $k\to\infty$ we get
$\mbf{P}\{\nu^h(2\varepsilon)< \infty\} = 1$. By repeating the above
procedure successively on the intervals $[h,2h]$, $[2h,3h], \ldots$ we
get $\mbf{P}\{\nu^a(2\varepsilon)< \infty\} = 1$ for every $a>0$. Let
$\Omega_1 = \bigcap_{m=1}^\infty\{\nu^m(1/m)<\infty\}$. Then
$\Omega_1\in\mcr{G}$ and \mbox{$\mbf{P}(\Omega_1)=1$}. Moreover, for
$\omega\in\Omega_1$ we can define a $\rho$-c\`{a}dl\`{a}g path
$t\mapsto
Y_t(\omega)$ on $[0,\infty)$ by $Y_t(\omega) := \lim_{F\ni s\to t+}
X_s(\omega)$. Take $x_0\in E$ and define $Y_t(\omega) = x_0$ for
$t\ge0$
and $\omega\in\Omega\setminus\Omega_1$. By (\ref{3.10}) one
can see
$t\mapsto X_t$ is right continuous in probability, so $Y_t = X_t$ a.s.
for every $t\ge0$. Then $\{Y_t\dvtx t\ge0\}$ is a $\rho$-c\`{a}dl\`{a}g
modification of $\{X_t\dvtx t\ge0\}$.
\end{pf}

Let $D[0,\infty)$ be the space of nonnegative c\`{a}dl\`{a}g functions
on $[0,\infty)$, and let $\mcr{B}(D[0,\infty))$ be its Borel
$\sigma$-algebra generated by the Skorokhod topology.
Theorems \ref{t3.2} and \ref{t3.3} imply that $\{Y(v)\dvtx v\ge0\}$ is
a nondecreasing process in $(D[0,\infty),
\mcr{B}(D[0,\infty)))$\vadjust{\goodbreak}
with independent increments. Let $\rho$ be the metric on
$D[0,\infty)$ defined by
%
%
\begin{equation}\label{3.11}
\rho(\xi,\zeta) = \int_0^\infty e^{-t}\sup_{0\le s\le t}
\bigl(|\xi(s)-\zeta(s)|\land1\bigr)\,dt.
\end{equation}
This metric corresponds to the topology of local uniform
convergence, which is strictly stronger than the Skorokhod topology.
\begin{theorem}\label{t3.6} The path-valued process $\{Y(v)\dvtx v\ge
0\}$ has
a $\rho$-c\`{a}dl\`{a}g~mo\-dification. Consequently, there is a version
of the solution flow $\{Y_t(v)\dvtx t\ge0,\allowbreak v\ge0\}$ of (\ref{3.5})
with the following properties:
\begin{longlist}
\item
for each $v\ge0$, $t\mapsto Y_t(v)$ is a c\`{a}dl\`{a}g process on
$[0,\infty)$ and solves~(\ref{3.5});

\item
for each $t\ge0$, $v\mapsto Y_t(v)$ is a nonnegative and
nondecreasing c\`{a}dl\`{a}g process on $[0,\infty)$.
\end{longlist}
\end{theorem}
\begin{pf}
\textit{Step} 1. For any $T\ge0$ let $D[0,T]$ be the space of
nonnegative c\`{a}dl\`{a}g functions on $[0,T]$, and let $\mcr
{B}(D[0,T])$ be
its $\sigma$-algebra generated by the Skorokhod topology. For $v\ge0$
let $Y^T(v) = \{Y_t(v)\dvtx 0\le t\le T\}$. Theorem~\ref{t3.3} implies that
$\{Y^T(v)\dvtx v\ge0\}$ is a process in $(D[0,T],\mcr{B}(D[0,T]))$ with
independent increments.

\textit{Step} 2. Let $F_T = \{T,r_1,r_2,\ldots\}$ be a countable
dense subset of $[0,T]$. We consider the metric $\rho_T$ on $D[0,T]$
defined by
\[
\rho_T(\xi,\zeta)
=
\sup_{0\le s\le T} |\xi(s)-\zeta(s)|
=
\sup_{r\in F_T} |\xi(s)-\zeta(s)|.
\]
For any $\varepsilon>0$ and $\xi\in D[0,T]$ we have
\begin{eqnarray*}
\bar{U}_\varepsilon(\xi)
:\!&=&
\{\zeta\in D[0,T]\dvtx \rho_T(\xi,\zeta)\le\varepsilon\} \\
&=&
\bigcap_{r\in F_T}\{\zeta\in D[0,T]\dvtx |\xi_r-\zeta_r|\le
\varepsilon\}.
\end{eqnarray*}
Then the above set belongs to $\mcr{B}(D[0,T])$ [see, e.g., Ethier
and Kurtz (\citeyear{EtK86}), pa\-ge~127]. It follows that
\[
U_\varepsilon(\xi)
:=
\{\zeta\in D[0,T]\dvtx \rho_T(\xi,\zeta)< \varepsilon\}
=
\bigcup_{n=1}^\infty\bar{U}_{\varepsilon-1/n}(\xi)
\]
also belongs to $\mcr{B}(D[0,T])$.

\textit{Step} 3. Let $(\mcr{F}_v^T)_{v\ge0}$ be the natural
filtration of $\{Y^T(v)\dvtx v\ge0\}$. For any $\varepsilon>0$ and
$0\le
s,t\le v$ we have
\[
\rho_T(Y^T(s),Y^T(t))
=
\sup_{r\in F_T} |Y_r(s)-Y_r(t)|.
\]
Then one can show $\{\omega\in\Omega\dvtx \rho_T(Y^T(\omega,s),
Y^T(\omega,t))< \varepsilon\}\in\mcr{F}_v^T$.

\textit{Step} 4. Let $(P_{u,v}^T)_{v\ge u}$ denote the transition
semigroup of $\{Y^T(v)\dvtx v\ge0\}$. By Proposition \ref{p3.4} for
$\varepsilon>0$ and $\xi\in D[0,\infty)$ we have
\begin{eqnarray*}
P_{u,v}(\xi,U_\varepsilon(\xi)^c)
&=&
\mbf{P}\Bigl\{\sup_{0\le s\le T}[Y_s(v) - Y_s(u)]\ge\varepsilon\Bigr\}
\\[-2pt]
&\le&
\varepsilon^{-1}\mbf{E}\Bigl\{\sup_{0\le s\le T}[Y_s(v) - Y_s(u)]\Bigr\}
\\[-2pt]
&\le&
\varepsilon^{-1}C(t)\bigl\{(v-u) + [\gamma(v)-\gamma(u)] \\
&&\hspace*{40.7pt}{}
+ \sqrt{v-u}+ \sqrt{\gamma(v)-\gamma(u)}\bigr\}.
\end{eqnarray*}
Since $v\mapsto\gamma(v)$ is uniformly continuous on each bounded
interval, Lemma~\ref{l3.5} implies that $\{Y^T(v)\dvtx v\ge0\}$ has a
$\rho_T$-c\`{a}dl\`{a}g modification. That implies the existence of a
$\rho$-c\`{a}dl\`{a}g modification of $\{Y(v)\dvtx v\ge0\}$.
\end{pf}

In the situation of Theorem \ref{t3.6} we call the solution
$\{Y_t(v)\dvtx t\ge0, v\ge0\}$ of (\ref{3.5}) a \textit{flow of
CBI-processes}. Let $F[0,\infty)$ be the set of nonnegative and
nondecreasing c\`{a}dl\`{a}g functions on $[0,\infty)$. Given a finite
stopping time $\tau$ and a function $\mu\in F[0,\infty)$ let
$\{Y^\mu_{\tau,t}(v)\dvtx t\ge0\}$ be the solution of
%
%
\begin{eqnarray}\label{3.12}
Y^\mu_{\tau,t}(v) &=&\mu(v) + \sigma\int_\tau^{\tau+t}
\int_0^{Y^\mu_{\tau,s-}(v)} W(ds,du) \nonumber\\[-2pt]
&&{}
+ \int_\tau^{\tau+t} [\gamma(v)-bY^\mu_{\tau,s-}(v)]\,ds \\[-2pt]
&&{}
+ \int_\tau^{\tau+t}\int_0^\infty\int_0^{Y^\mu_{\tau,s-}(v)} z
\tilde{N}(ds,dz,du)\nonumber
\end{eqnarray}
and write simply $\{Y^\mu_t(v)\dvtx t\ge0\}$ instead of $\{Y^\mu
_{0,t}(v)\dvtx
t\ge0\}$. The pathwise uniqueness for the above equation follows from
that of (\ref{3.5}) since $\{W(\tau+ds,du)\}$ is a white noise based on
$ds\,dz$, and $\{N(\tau+ds,dz,du)\}$ is a Poisson random measure with
intensity $ds\,m(dz)\,du$. Let $G_{\tau,t}$ be the random operator on
$F[0,\infty)$ that maps $\mu$ to $Y^\mu_{\tau,t}$.
\begin{theorem}\label{t3.7} For any finite stopping time $\tau$ we have
$\mbf{P}\{Y^\mu_{\tau+t} = G_{\tau,t}Y^\mu_\tau$ for all $t\ge
0\} = 1$.
\end{theorem}
\begin{pf}
By the sample path regularity of $(t,v)\mapsto Y_t(v)$ we only
need to show $\mbf{P}\{Y^\mu_{\tau+t}(v) = G_{\tau,t}Y^\mu_\tau
(v)\} = 1$
for every $t\ge0$ and $v\ge0$. By (\ref{3.5}) we have
\begin{eqnarray*}
Y^\mu_{\tau+t}(v) &=& Y^\mu_\tau(v) + \sigma\int_\tau^{\tau+t}
\int_0^{Y^\mu_{s-}(v)} W(ds,du) \\[-2pt]
&&{}
+ \int_\tau^{\tau+t}[\gamma(v)-bY^\mu_{s-}(v)]\,ds \\[-2pt]
&&{}
+ \int_\tau^{\tau+t} \int_0^\infty\int_0^{Y^\mu_{s-}(v)} z
\tilde{N}(ds,dz,du).
\end{eqnarray*}
By the pathwise uniqueness for (\ref{3.12}) we get the desired result.
\end{pf}

For any Radon measure $\mu(dv)$ on $[0,\infty)$ with distribution
function \mbox{$v\,{\mapsto}\,\mu(v)$}, the random function $v\mapsto
Y^\mu_t(v)$ induces a random Radon measu\-re~$Y^\mu_t(dv)$ on
$[0,\infty)$ so that $Y^\mu_t([0,v]) = Y^\mu_t(v)$ for $v\ge0$. We
shall give some characterizations of the measure-valued process
$\{Y^\mu_t\dvtx t\ge0\}$.

For simplicity, we fix a constant $a\ge0$ and consider the restrictions
of~$\mu(dv)$, $\gamma(dv)$ and $\{Y^\mu_t\dvtx t\ge0\}$ to $[0,a]$ without
changing the notation. Let us consider the step function
%
%
\begin{equation}\label{3.13}
f(x) = c_01_{\{0\}}(x) + \sum_{i=1}^n c_i1_{(a_{i-1},a_i]}(x),\qquad
x\in[0,a],
\end{equation}
where $\{c_0,c_1,\ldots,c_n\}\subset\mbb{R}$ and $\{0=a_0 <a_1<
\cdots<
a_n=a\}$ is a partition of $[0,a]$. For this function we have
%
%
\begin{equation}\label{3.14}
\langle Y^\mu_t,f\rangle= c_0Y^\mu_t(0) + \sum_{i=1}^n c_i[Y^\mu
_t(a_i) -
Y^\mu_t(a_{i-1})].
\end{equation}
{From} (\ref{3.12}) and (\ref{3.14}) it is simple to see
%
%
\begin{eqnarray}\label{3.15}
\langle Y^\mu_t,f\rangle
&=&
\langle\mu,f\rangle+ \sigma\int_0^t\int_0^\infty g^\mu
_{s-}(u)W(ds,du)\nonumber\\
&&{}
+ \int_0^t [\langle\gamma,f\rangle- b\langle Y^\mu
_{s-},f\rangle] \,ds \\
&&{}
+ \int_0^t\int_0^\infty\int_0^\infty zg^\mu_{s-}(u)
\tilde{N}(ds,dz,du),\nonumber
\end{eqnarray}
where
%
%
\begin{equation}\label{3.16}
g^\mu_s(u) = c_01_{\{u\le Y^\mu_s(0)\}} + \sum_{i=1}^n
c_i1_{\{Y^\mu_s(a_{i-1})< u\le Y^\mu_s(a_i)\}}.
\end{equation}
\begin{prop}\label{p3.8} For any $t\ge0$ and $f\in B[0,a]$ we have
%
%
\begin{equation}\label{3.17}
\mbf{E}[\langle Y^\mu_t,f\rangle] = \langle\mu,f\rangle e^{-bt} +
\langle\gamma,f\rangle b^{-1}(1 -
e^{-bt})
\end{equation}
with $b^{-1}(1 - e^{-bt}) = t$ for $b=0$ by convention.
\end{prop}
\begin{pf}
We first consider the step function (\ref{3.13}). By taking the
expectation in (\ref{3.15}) we obtain
\[
\mbf{E}[\langle Y^\mu_t,f\rangle]
=
\langle\mu,f\rangle+ t\langle\gamma,f\rangle- b\int_0^t \mbf
{E}[\langle Y^\mu_s,f\rangle] \,ds.\vadjust{\goodbreak}
\]
The above integral equation has the unique solution given by
(\ref{3.17}). For a~general function $f\in B[0,a]$ we get (\ref
{3.17}) by
a monotone class argument.
\end{pf}
\begin{theorem}\label{t3.9} The measure-valued process $\{Y^\mu_t\dvtx
t\ge0\}$ is
a c\`{a}dl\`{a}g strong Markov process in $M[0,a]$ with $Y^\mu_0 = \mu$.
\end{theorem}
\begin{pf}
In view of (\ref{3.14}), the process $t\mapsto\langle Y^\mu
_t,f\rangle$ is
c\`{a}dl\`{a}g for the step function (\ref{3.13}). Since any function in
$C[0,a]$ can be approximated by a~sequence of step functions in the
supremum norm, it is easy to conclude $t\mapsto\langle Y^\mu
_t,f\rangle$ is
c\`{a}dl\`{a}g for all $f\in C[0,a]$. By Theorem \ref{t3.7}, for any finite
stopping time $\tau$ we have $Y^\mu_{\tau+t} = G_{\tau,t} Y^\mu
_\tau$
almost surely.\vspace*{1pt} That clearly implies the strong Markov property of
$\{Y^\mu_t\dvtx t\ge0\}$.
\end{pf}
\begin{theorem}\label{t3.10} For any $f\in B[0,a]$ the process
$\{\langle Y^\mu_t,f\rangle\dvtx t\ge0\}$ has a~c\`{a}dl\`{a}g
modification. Moreover,
there is a locally bounded function $t\mapsto C(t)$ so that
%
%
\begin{equation}\label{3.18}\qquad
\mbf{E}\Bigl[\sup_{0\le s\le t}\langle Y^\mu_s,f\rangle\Bigr]
\le
C(t)[\langle\mu,f\rangle+ \langle\gamma,f\rangle
+ \langle\mu,f^2\rangle^{1/2} + \langle\gamma,f^2\rangle
^{1/2}]
\end{equation}
for every $t\ge0$ and $f\in B[0,a]^+$.
\end{theorem}
\begin{pf}
We first consider a nonnegative step function given by
(\ref{3.13}) with constants $\{c_0,c_1,\ldots,c_n\}\subset
\mbb{R}_+$. By (\ref{3.15}) and Doob's martingale inequality,
\begin{eqnarray*}
&&\mbf{E}\Bigl[\sup_{0\le s\le t}\langle Y^\mu_s,f\rangle\Bigr]\\
&&\qquad\le
\langle\mu,f\rangle+ 2\sigma\mbf{E}^{1/2} \biggl\{
\biggl[\int_0^t
\int_0^\infty g^\mu_{s-}(u) W(ds,du)\biggr]^2\biggr\} \\
&&\qquad\quad{}
+ t\langle\gamma,f\rangle+ |b|\int_0^t\mbf{E}[\langle Y_s^\mu
,f\rangle]\,ds \\
&&\qquad\quad{}
+ 2\mbf{E}^{1/2}\biggl\{\biggl[\int_0^t \int_0^1 \int
_0^\infty z
g^\mu_{s-}(u) \tilde{N}(ds,dz,du)\biggr]^2\biggr\} \\
&&\qquad\quad{}
+ \mbf{E}\biggl[\int_0^t \int_1^\infty\int_0^\infty z g^\mu_{s-}(u)
N(ds,dz,du)\biggr] \\
&&\qquad=
\langle\mu,f\rangle+ 2\sigma\mbf{E}^{1/2} \biggl[\int_0^tds
\int_0^\infty g^\mu_s(u)^2 \,du\biggr] \\
&&\qquad\quad{}
+ t\langle\gamma,f\rangle+ |b|\int_0^t\mbf{E}[\langle Y_s^\mu
,f\rangle]\,ds \\
&&\qquad\quad{}
+ 2\mbf{E}^{1/2}\biggl[\int_0^tds \int_0^1 z^2 m(dz)
\int_0^\infty g^\mu_s(u)^2 \,du\biggr] \\
&&\qquad\quad{}
+ \mbf{E}\biggl[\int_0^tds \int_1^\infty z m(dz) \int_0^\infty
g^\mu_s(u) \,du\biggr] \\
&&\qquad\le
\langle\mu,f\rangle+ 2\biggl(\int_0^t\mbf{E}[\langle Y_s^\mu
,f^2\rangle]\,ds\biggr)^{1/2}
\biggl[\sigma+ \biggl(\int_0^1 z^2 m(dz)\biggr)^{1/2}\biggr]
\\
&&\qquad\quad{}
+ t\langle\gamma,f\rangle+ \int_0^t\mbf{E}[\langle Y_s^\mu
,f\rangle]\,ds\,\biggl[|b| +
\int_1^\infty z m(dz)\biggr].
\end{eqnarray*}
In view of (\ref{3.17}) we get (\ref{3.18}) for the step function.
Now let $\eta(dv) = \mu(dv) + \gamma(dv)$ and choose a bounded
sequence of step functions $\{f_n\}$ so that $f_n\to f$ in
$L^2(\eta)$ as $n\to\infty$. By applying (\ref{3.18}) to the
nonnegative step function $|f_n-f_m|$ we get
\[
\mbf{E}\Bigl[\sup_{0\le s\le t}\langle Y^\mu_s,|f_n-f_m|\rangle\Bigr]
\le
C(t)[\langle\eta,|f_n-f_m|\rangle+ 2\langle\eta
,|f_n-f_m|^2\rangle^{1/2}].
\]
The right-hand side tends to zero as $m,n\to\infty$. Then there is
a c\`{a}dl\`{a}g process $\{Y^\mu_t(f)\dvtx t\ge0\}$ so that
%
%
\begin{equation}\label{3.19}
\mbf{E}\Bigl[\sup_{0\le s\le t}|\langle Y^\mu_s,f_n\rangle-
Y^\mu_s(f)|\Bigr]\to0,\qquad n\to\infty.
\end{equation}
On the other hand, from (\ref{3.17}) we have
\[
\mbf{E}[\langle Y^\mu_t,|f_n-f|\rangle]
=
\langle\mu,|f_n-f|\rangle e^{-bt} + b^{-1}(1 - e^{-bt})\langle\gamma
,|f_n-f|\rangle,
\]
which tends to zero as $n\to\infty$. Then $\{Y^\mu_t(f)\dvtx t\ge0\}$
is a modification of $\{\langle Y^\mu_t,f\rangle\dvtx t\ge0\}$. Finally,
we get
(\ref{3.18}) for $f\in B[0,a]^+$ by using (\ref{3.19}) and the
result for step functions.
\end{pf}
\begin{theorem}\label{t3.11} The process $\{Y^\mu_t\dvtx t\ge0\}$ is the
unique solution of the following martingale problem: for every $G\in
C^2(\mbb{R})$ and $f\in B[0,a]$,
%
%
\begin{eqnarray}\label{3.20}
&&G(\langle Y^\mu_t,f\rangle)\nonumber\\
&&\qquad=
G(\langle\mu,f\rangle) + \frac{1}{2}\sigma^2 \int_0^t
G^{\prime\prime}(\langle Y^\mu_s,f\rangle) \langle Y^\mu
_s,f^2\rangle\,ds\nonumber\\
&&\qquad\quad{}
+ \int_0^t G^\prime(\langle Y^\mu_s,f\rangle)[\langle\gamma
,f\rangle- b\langle Y^\mu_s,f\rangle] \,ds
\nonumber\\[-8pt]\\[-8pt]
&&\qquad\quad{}
+ \int_0^t ds\int_{[0,a]} Y^\mu_s(dx)
\int_0^\infty\bigl[G\bigl(\langle Y^\mu_s,f\rangle+ zf(x)\bigr)\nonumber\\
&&\qquad\quad\hspace*{120.3pt}{}
- G(\langle Y^\mu_s,f\rangle) - zf(x)G^\prime(\langle Y^\mu
_s,f\rangle)\bigr]m(dz) \hspace*{-14pt}\nonumber\\
&&\qquad\quad{}
+ \mbox{local mart.}\nonumber
\end{eqnarray}
\end{theorem}
\begin{pf}
Again we start with the step function (\ref{3.13}). Using
(\ref{3.15}) and It\^{o}'s formula,
\begin{eqnarray*}
&&G(\langle Y^\mu_t,f\rangle)\\
&&\qquad=
G(\langle\mu,f\rangle) + \frac{1}{2}\sigma^2 \int_0^tds\int
_0^\infty
G^{\prime\prime}(\langle Y^\mu_{s-},f\rangle) g^\mu_{s-}(u)^2 \,du
\\
&&\qquad\quad{}
+ \int_0^t G^\prime(\langle Y^\mu_{s-},f\rangle)[\langle\gamma
,f\rangle- b
\langle Y^\mu_{s-},f\rangle]\,ds \\
&&\qquad\quad{}
+ \int_0^tds\int_0^\infty m(dz)\int_0^\infty\bigl[G\bigl(\langle Y^\mu
_s,f\rangle+
zg^\mu_s(u)\bigr) \\
&&\qquad\quad\hspace*{111.3pt}{}
- G(\langle Y^\mu_s,f\rangle) - G^\prime(\langle Y^\mu_s,f\rangle
) zg^\mu_s(u)\bigr]\,du \\
&&\qquad\quad{}+
\mbox{local mart.} \nonumber\\
&&\qquad=
G(\langle\mu,f\rangle) + \frac{1}{2}\sigma^2 \int_0^t
G^{\prime\prime}(\langle Y^\mu_s,f\rangle) \langle Y^\mu
_s,f^2\rangle\,ds \\
&&\qquad\quad{}
+ \int_0^t G^\prime(\langle Y^\mu_s,f\rangle)[\langle\gamma
,f\rangle- b \langle Y^\mu_s,f\rangle]\,ds
\\
&&\qquad\quad{}
+ \int_0^tds\int_0^\infty Y^\mu_s(dx)\int_0^\infty
\bigl[G\bigl(\langle Y^\mu_s,f\rangle+ zf(x)\bigr) \\
&&\qquad\quad\hspace*{118.2pt}{}
- G(\langle Y^\mu_s,f\rangle) - G^\prime(\langle Y^\mu_s,f\rangle
) zf(x)\bigr]m(dz) \\
&&\qquad\quad{}+
\mbox{local mart.}
\end{eqnarray*}
That proves (\ref{3.20}) for step functions. For $f\in B[0,a]$ we get the
martingale problem using (\ref{3.19}). The uniqueness of the solution
follows from a result in Li (\citeyear{Li10}), Section 9.3.
\end{pf}

The solution of the martingale problem (\ref{3.20}) is the special
case of
the \textit{immigration superprocess} studied in \citet{Li10}
with trivial
spatial motion. More precisely, the infinitesimal particles propagate in
$[0,a]$ without migration. Then for any disjoint bounded Borel subsets
$B_1$ and $B_2$ of $[0,a]$, the nonnegative real-valued processes
$\{Y^\mu_t(B_1)\dvtx t\ge0\}$ and $\{Y^\mu_t(B_2)\dvtx t\ge0\}$ are
independent.
That explains why the restriction of $\{Y^\mu_t\dvtx t\ge0\}$ to the interval
$[0,a]$ is still a Markov process. To consider the process of measures on the half
line $[0,\infty)$ we need to introduce a weight function as follows.

Let $h$ be a strictly positive continuous function on $[0,\infty)$
vanishing at infinity. Let $M_h[0,\infty)$ be the space of Radon measures
$\mu$ on $[0,\infty)$ so that $\langle\mu,h\rangle<\infty$. Let
$B_h[0,\infty)$ be
the set of Borel functions on $[0,\infty)$ bounded by $\const\cdot
h$, and
let $C_h[0,\infty)$ denote its subset of continuous functions. A~topo\-logy
on $M_h[0,\infty)$ can be defined by the convention $\mu_n\to\mu$ in
$M_h[0,\infty)$ if and only\vadjust{\goodbreak} if $\langle\mu_n,f\rangle\to\langle
\mu,f\rangle$ for every $f\in
C_h[0,\infty)$. Suppose that $\mu\in M_h[0,\infty)$ and $\gamma\in
M_h[0,\infty)$. It is easy to show that $\{Y^\mu_t\dvtx t\ge0\}$ is a
c\`{a}dl\`{a}g strong Markov process in $M_h[0,\infty)$, and the
results of
Theorem \ref{t3.10} and Theorem~\ref{t3.11} are also true for
$B_h[0,\infty)$.


\section{Generalized Fleming--Viot flows}\label{sec4}

In this section we give a construction of the generalized Fleming--Viot
flow as the strong solution of a stochastic integral equation. Let
$\sigma\ge0$, $b\ge0$ and $0\le\beta\le1$ be constants, and
let~$z^2\nu(dz)$ be a finite measure on $(0,1]$. Suppose that $\{B(ds,du)\}$
is a~white noise on $(0,\infty)^2$ with intensity $ds\,du$, and
$\{M(ds,dz,du)\}$ is a Poisson random measure on $(0,\infty)\times
(0,1]\times(0,\infty)$ with intensity $ds\,\nu(dz)\,du$. Let
\[
q(x,u) = 1_{\{u\le1\land x\}} - (1\land x),\qquad x\ge0, u\in
(0,1].
\]
We first consider the stochastic integral equation
%
%
\begin{eqnarray}\label{4.1}
X_t &=& X_0 + \int_0^t\int_0^1\sigma q(X_{s-},u) B(ds,du) +
\int_0^t b(\beta-X_{s-}) \,ds \nonumber\\[-8pt]\\[-8pt]
&&{}
+ \int_0^t\int_0^1\int_0^1 zq(X_{s-},u) \tilde{M}(ds,dz,du),\nonumber
\end{eqnarray}
where $\tilde{M}(ds,dz,du)$ denotes the compensated measure of
$M(ds,dz,du)$. In fact, the compensation in (\ref{4.1}) can be
disregarded as
\[
\int_0^1 q(X_{s-},u)\,du = \int_0^1 \bigl[1_{\{u\le X_{s-}\land1\}} -
(X_{s-}\land1)\bigr]\,du = 0.
\]
\begin{theorem}\label{t4.1} There is a unique nonnegative strong solution
to (\ref{4.1}).
\end{theorem}
\begin{pf}
We first show the pathwise uniqueness for (\ref{4.1}). Set
$l(x,y,u) = q(x,u) - q(y,u)$. For $x,y\ge0$ and $0\le z,t\le1$ we have
\begin{eqnarray*}
&&(x-y) + ztl(x,y,u) \nonumber\\
&&\qquad
= [(x-1\land x) - (y-1\land y)] + (1-zt)(1\land x-1\land y) \nonumber\\
&&\qquad\quad{}
+ zt\bigl(1_{\{u\le x\land1\}}-1_{\{u\le y\land1\}}\bigr).
\end{eqnarray*}
It is then easy to see
\[
|(x-y) + ztl(x,y,u)|
\ge
(1-zt)|1\land x-1\land y|.
\]
Moreover, we have
\begin{eqnarray*}
\int_0^1 l(x,y,u)^2\,du
&=&
(1\land x-1\land y) - (1\land x-1\land y)^2 \\
&\le&
|1\land x-1\land y|.
\end{eqnarray*}
Using the above two inequalities,
\begin{eqnarray*}
&&\int_0^1(1-t)\,dt\int_0^1\nu(dz)\int_0^1 \frac{z^2l(x,y,u)^2}
{|(x-y) + ztl(x,y,u)|}\,du \\
&&\qquad
\le\int_0^1z^2\nu(dz) \int_0^1\frac{1-t}{1-zt}\,dt\int_0^1
\frac{l(x,y,u)^2} {|1\land x-1\land y|}\,du \\
&&\qquad
\le\int_0^1z^2\nu(dz) \int_0^1\frac{1-t}{1-zt}\,dt \le\int_0^1z^2
\nu(dz).
\end{eqnarray*}
Then condition \textup{(2.d)} is satisfied with $\rho(z) = \sqrt{z}$. Other
conditions of Theorem \ref{t2.1} can be checked easily. Then we have the
pathwise uniqueness for~(\ref{4.1}). To show the existence of the
solution, we may assume $X_0 = v\ge0$ is a~deterministic constant. By
Theorem \ref{t2.5} there a unique nonnegative strong solution of
(\ref{4.1}) if the Poisson integral term is removed. Then for each
$k\ge
1$ there is a unique nonnegative strong solution to
%
%
\begin{eqnarray}\label{4.2}
Z_t &=& Z_0 + \int_0^t\int_0^1 \sigma q(Z_{s-},u) B(ds,du)\nonumber\\
&&{} +
\int_0^t b(\beta-Z_{s-}) \,ds \\
&&{}
+ \int_0^t\int_{1/k}^1\int_0^1 zq(Z_{s-},u) M(ds,dz,du),\nonumber
\end{eqnarray}
because the last term on the right-hand side gives at most a finite
number of jumps on each bounded time interval. Let $\{Z_k(t)\dvtx t\ge
0\}$
be the solution of (\ref{4.2}) with $Z_k(0) = v$. Let $T_1 = \inf\{
t\ge
0\dvtx Z_k(t)\le1\}$. On the time interval $[0,T_1]$, the stochastic
integral terms in (\ref{4.2}) vanish. Then $t\mapsto Z_k(t)$ is
nonincreasing on $[0,T_1]$. By modifying the proof of Proposition 2.1 in
\citet{FuL10} one can see $Z_k(t)\le1$ for $t\ge T_1$. Thus
$Z_k(t)\le(Z_k(0)\vee1) = (v\vee1)$ for all $t\ge0$. Let $\{\tau
_k\}$
be a bounded sequence of stopping times. Note that the last term on the
right-hand side of (\ref{4.2}) can be considered as a stochastic integral
with respect to the compensated Poisson random measure. Then for any
$t\ge0$ we have
\begin{eqnarray*}
&&\mbf{E}\{[Z_k(\tau_k+t)-Z_k(\tau_k)]^2\} \\
&&\qquad
\le3\sigma^2\mbf{E}\biggl[\int_0^tds\int_0^1q\bigl(Z_k(\tau_k+s),u\bigr)^2\,
du\biggr]
+ 3b^2t^2(v\vee1)^2 \\
&&\qquad\quad{}
+ 3\mbf{E}\biggl[\int_0^tds\int_0^1z^2\nu(dz)\int_0^1
q\bigl(Z_k(\tau_k+s),u\bigr)^2 \,du\biggr] \\
&&\qquad
\le3t\biggl[\sigma^2 + tb^2(v\vee1)^2 + \int_0^1z^2\nu(dz)\biggr].
\end{eqnarray*}
The right-hand side tends to zero as $t\to0$. By a criterion of
\citet{Ald78}, the sequence $\{Z_k(t)\dvtx t\ge0\}$ is tight in
$D([0,\infty),\mbb{R}_+)$ [see also Ethier and Kurtz
(\citeyear{EtK86}), pages 137 and 138].
By a modification of the proof of Theorem~4.4 in \citet{FuL10}
one sees that any limit point of this sequence is a~weak
solution of (\ref{4.1}).\vspace*{-3pt}
\end{pf}

Now let $v\mapsto\gamma(v)$ be a nondecreasing continuous function
on $[0,1]$ so that $0\le\gamma(v)\le1$ for all $0\le v\le1$. We
denote by $\gamma(dv)$ the sub-probability measure on $[0,1]$ so
that $\gamma([0,v]) = \gamma(v)$ for $0\le v\le1$. By
Theorem \ref{t4.1} for each $v\ge0$ there is a pathwise unique
nonnegative solution $\{X_t(v)\dvtx t\ge0\}$ to the equation
%
%
\begin{eqnarray}\label{4.3}
X_t(v) &=& v + \int_0^t\int_0^1 \sigma\bigl[1_{\{u\le X_{s-}(v)\}} -
X_{s-}(v)\bigr] B(ds,du) \nonumber\\[-2pt]
&&{}
+ \int_0^t b[\gamma(v)-X_{s-}(v)] \,ds \\[-2pt]
&&{}
+ \int_0^t\int_0^1\int_0^1z\bigl[1_{\{u\le X_{s-}(v)\}} - X_{s-}(v)\bigr]
\tilde{M}(ds,dz,du).\nonumber
\end{eqnarray}
It is not hard to see that $0\le v\le1$ implies $\mbf{P}\{0\le
X_t(v)\le
1$ for all $t\ge0\} = 1$. The compensation for the Poisson random
measure can be disregarded, so this equation just coincides with
(\ref{1.6}). By Theorem \ref{t2.2} for any $0\le v_1\le v_2\le1$ we have
\[
\mbf{P}\{X_t(v_1)\le X_t(v_2) \mbox{ for all } t\ge0\} = 1.
\]
Therefore $\{X(v)\dvtx 0\le v\le1\}$ is a nondecreasing path-valued process
in $D[0,\infty)$.

\begin{prop}\label{p4.2} There is a locally bounded nonnegative
function $t\mapsto C(t)$ on $[0,\infty)$ so that
%
%
\begin{eqnarray}\label{4.4}
&&
\mbf{E}\Bigl\{\sup_{0\le s\le t}[X_s(v_2) - X_s(v_1)]\Bigr\}\nonumber\\[-2pt]
&&\qquad\le
C(t)\bigl\{(v_2-v_1) + [\gamma(v_2)-\gamma(v_1)] \\[-2pt]
&&\qquad\quad\hspace*{23.5pt}{}
+ \sqrt{v_2-v_1} + \sqrt{\gamma(v_2)-\gamma(v_1)}\bigr\}\nonumber
\end{eqnarray}
for $t\ge0$ and $0\le v_1\le v_2\le1$.
\end{prop}
\begin{pf}
Let $Z_t = X_t(v_2)-X_t(v_1)$ for $t\ge0$. {From}
(\ref{4.3}) we have
%
%
\begin{eqnarray}\label{4.5}
Z_t &=&(v_2-v_1) + \int_0^t\int_0^1 \sigma[Y_{s-}(u)-Z_{s-}]
B(ds,du) \nonumber\\[-2pt]
&&{}
+ \int_0^t b\{[\gamma(v_2)-\gamma(v_1)] - Z_{s-}\} \,ds \\[-2pt]
&&{}
+ \int_0^t\int_0^1\int_0^1 z[Y_{s-}(u)-Z_{s-}]
\tilde{M}(ds,dz,du),\nonumber
\end{eqnarray}
where $Y_s(u) = 1_{\{X_s(v_1)< u\le X_s(v_2)\}}$. Taking the expectation
in (\ref{4.5}) and solving a deterministic integral equation one can show
%
%
\begin{equation}\label{4.6}
\mbf{E}[Z_t] = (v_2-v_1)e^{-bt} + [\gamma(v_2)-\gamma(v_1)]
(1-e^{-bt}).
\end{equation}
By (\ref{4.5}) and Doob's martingale inequality,
\begin{eqnarray*}
\mbf{E}\Bigl\{\sup_{0\le s\le t}Z_s\Bigr\}
&\le&
(v_2-v_1) + 2\sigma\mbf{E}^{1/2}\biggl\{\biggl(\int_0^t
\int_0^1
[Y_{s-}(u)-Z_{s-}] B(ds,du)\biggr)^2\biggr\} \\[-2pt]
&&{}
+ \int_0^t b\{[\gamma(v_2)-\gamma(v_1)] + \mbf{E}[Z_s]\} \,ds \\[-2pt]
&&{}
+ 2\mbf{E}^{1/2}\biggl\{\biggl(\int_0^t\int_0^1\int_0^1
z[Y_{s-}(u)-Z_{s-}] \tilde{M}(ds,dz,du)\biggr)^2\biggr\} \\[-2pt]
&=&
(v_2-v_1) + 2\sigma\mbf{E}^{1/2}\biggl\{\int_0^tds \int_0^1
[Y_s(u) - Z_s]^2 \,du\biggr\} \\[-2pt]
&&{}
+ \int_0^t b\{[\gamma(v_2)-\gamma(v_1)] + \mbf{E}[Z_s]\} \,ds \\[-2pt]
&&{}
+ 2\mbf{E}^{1/2}\biggl\{\int_0^tds\int_0^1 z^2\nu(dz)
\int_0^1 [Y_s(u) - Z_s]^2 \,du\biggr\},
\end{eqnarray*}
where
\[
\int_0^1 [Y_s(u) - Z_s]^2 \,du
=
Z_s(1 - Z_s) \le Z_s.
\]
Then we have (\ref{4.4}) by (\ref{4.6}).
\end{pf}

Recall that $D[0,\infty)$ is the space of nonnegative c\`{a}dl\`{a}g
functions on $[0,\infty)$ endowed with the Borel $\sigma$-algebra
generated by the Skorokhod topology. Let~$\rho$ be the metric on
$D[0,\infty)$ defined by (\ref{3.11}).
\begin{theorem}\label{t4.3} The path-valued process $\{X(v)\dvtx 0\le
v\le
1\}$ is a Markov process in $D[0,\infty)$.
\end{theorem}
\begin{pf}
Let $0<v<1$, and let $\tau_n = \inf\{t\ge0\dvtx X_t(v)\le1/n\}
$ for
$n\ge1$. In view of (\ref{4.3}), we have $X_t(v)=0$ if $X_{t-}(v)=0$.
Then $\tau_n\to\tau_\infty:= \inf\{t\ge0\dvtx\break X_t(v)=0\}$ as $n\to
\infty$. For any $p\in[0,v)$ the comparison property and pathwise
uniqueness for (\ref{4.3}) imply $X_t(p) = X_t(v)$ for $t\ge
\tau_\infty$. Let $Z_n(t) = X_{t\land\tau_n}(v)^{-1} X_{t\land
\tau_n}(p)$ for $t\ge0$. By (\ref{4.3}) and It\^{o}'s formula,
\begin{eqnarray*}
Z_n(t) &=&\frac{p}{v} + \int_0^{t\land\tau_n}\int_0^1
\frac{\sigma}{X_{s-}(v)}\bigl[1_{\{u\le X_{s-}(p)\}} - X_{s-}(p)\bigr]
B(ds,du) \\[-2pt]
&&{}
- \int_0^{t\land\tau_n}\int_0^1 \frac{\sigma X_{s-}(p)}
{X_{s-}(v)^2} \bigl[1_{\{u\le X_{s-}(v)\}} - X_{s-}(v)\bigr] B(ds,du)
\\[-2pt]
&&{}
+ \int_0^{t\land\tau_n} bX_{s-}(v)^{-1} [\gamma(p)
-\gamma(v)X_{s-}(v)^{-1}X_{s-}(p)] \,ds \\[-2pt]
&&{}
+ \int_0^{t\land\tau_n}ds\int_0^1\frac{\sigma^2X_{s-}(p)}
{X_{s-}(v)^3} \bigl[1_{\{u\le X_{s-}(v)\}} - X_{s-}(v)\bigr]^2\,du \\[-2pt]
&&{}
- \int_0^{t\land\tau_n}ds\int_0^1 \frac{\sigma^2}{X_{s-}(v)^2}
\bigl[1_{\{u\le X_{s-}(p)\}} - X_{s-}(p)\bigr] \\[-2pt]
&&\hspace*{68.2pt}{}
\times\bigl[1_{\{u\le X_{s-}(v)\}} - X_{s-}(v)\bigr]\, du \\[-2pt]
&&{}
+ \int_0^{t\land\tau_n}\int_0^1\int_0^1\biggl\{\frac{X_{s-}(p)(1-z)
+ z1_{\{u\le X_{s-}(p)\}}}{X_{s-}(v)(1-z) + z1_{\{u\le X_{s-}(v)\}}}
- \frac{X_{s-}(p)}{X_{s-}(v)}\biggr\}\\[-2pt]
&&\hspace*{70pt}{}\times M(ds,dz,du) \\[-2pt]
&=&
\frac{p}{v} + \int_0^{t\land\tau_n}\int_0^{X_{s-}(v)} \sigma
X_{s-}(v)^{-1}\bigl[1_{\{u\le X_{s-}(p)\}}
- X_{s-}(v)^{-1}X_{s-}(p)\bigr]\\[-2pt]
&&\hspace*{85.2pt}{} \times B(ds,du) \\[-2pt]
&&{}
+ \int_0^{t\land\tau_n} bX_{s-}(v)^{-1} [\gamma(p)
-\gamma(v)X_{s-}(v)^{-1}X_{s-}(p)] \,ds \\[-2pt]
&&{}
+ \int_0^{t\land\tau_n}\int_0^1\int_0^{X_{s-}(v)}
\biggl[\frac{X_{s-}(p)(1-z) + z1_{\{u\le X_{s-}(p)\}}}{X_{s-}(v)(1-z)
+ z}
- \frac{X_{s-}(p)}{X_{s-}(v)}\biggr] \\[-2pt]
&&\hspace*{91.5pt}{}\times M(ds,dz,du),
\end{eqnarray*}
where the two terms involving $\sigma^2$ counteract each other. Observe
also that the last integral does not change if we replace $M(ds,dz,du)$
by the compensated measure $\tilde{M}(ds,dz,du)$. Then we get the
equation
%
%
\begin{eqnarray}\label{4.7}\quad
Z_n(t) &=&\frac{p}{v} + \int_0^{t\land\tau_n}\int_0^{X_{s-}(v)}
\sigma X_{s-}(v)^{-1} \bigl[1_{\{u\le X_{s-}(v)Z_n(s-)\}}
- Z_n(s-)\bigr]\nonumber\\
&&\hspace*{85.3pt}{}\times B(ds,du) \nonumber\\
&&{}
+ \int_0^{t\land\tau_n}\int_0^1\int_0^{X_{s-}(v)}
z\biggl[\frac{1_{\{u\le X_{s-}(v)Z_n(s-)\}}}{z + (1-z)X_{s-}(v)}
- \frac{Z_n(s-)}{z + (1-z)X_{s-}(v)}\biggr] \\
&&\hspace*{92pt}{}\times\tilde{M}(ds,dz,du) \nonumber\\
&&{}
+ \int_0^{t\land\tau_n} bX_{s-}(v)^{-1}[\gamma(p) -\gamma
(v)Z_n(s-)] \,ds.\nonumber
\end{eqnarray}
Since $X_{s-}(v)\ge1/n$ for $0< s\le\tau_n$, by a simple generalization
of Theorem~\ref{t2.1} one can show the pathwise uniqueness holds for
(\ref{4.7}). Then, setting $Z_t = \lim_{n\to\infty}Z_n(t)$ we have
%
%
\begin{equation}\label{4.8}
X_t(p) = Z_tX_t(v)1_{\{t< \tau_\infty\}} + X_t(v)1_{\{t\ge
\tau_\infty\}},\qquad t\ge0.
\end{equation}
Now from (\ref{4.7}) and (\ref{4.8}) we infer that $\{X_t(p)\dvtx t\ge
0\}$
is measurable with respect to the $\sigma$-algebra $\mcr{F}_v$ generated
by the process $\{X_t(v)\dvtx t\ge0\}$ and the restricted martingale
measures
\[
1_{\{u\le X_{s-}(v)\}}B(ds,du),\qquad 1_{\{u\le X_{s-}(v)\}}
\tilde{M}(ds,dz,du).
\]
By similar arguments, for any $q\in(v,1]$ one can see $\{1-X_t(q)\dvtx
t\ge0\}$ is~mea\-surable with respect to the $\sigma$-algebra
$\mcr{G}_v$ generated by the process $\{1-X_t(v)\dvtx\allowbreak t\ge0\}$ and the
restricted martingale measures
\[
1_{\{X_{s-}(v)<u\le1\}}B(ds,du),\qquad 1_{\{X_{s-}(v)<u\le1\}}
\tilde{M}(ds,dz,du).
\]
Observe that $\{B(ds,X_{s-}(v)+du)\}$ is a white noise with intensity
$ds\,du$ and $\{M(ds,dz,X_{s-}(v)+du)\}$ is a Poisson random measure
with intensity $ds\,\hspace*{-0.5pt}\nu(dz)\hspace*{-0.5pt}\,du$. Then,
given $\{X_t(v)\dvtx t\ge 0\}$ the $\sigma$-algebras $\mcr{F}_v$ and
$\mcr{G}_v$ are conditionally independent. That implies the Markov
property of $\{(X(v),\mcr{F}_v)\dvtx\allowbreak 0\le v\le1\}$.
\end{pf}
\begin{theorem}\label{t4.4}
The path-valued Markov process $\{X(v)\dvtx 0\le v\le1\}$ has a
$\rho$-c\`{a}dl\`{a}g modification. Consequently, there is a version of
the solution flow $\{X_t(v)\dvtx t\ge0, 0\le v\le1\}$ of (\ref{4.3})
with the following properties:
\begin{longlist}
\item
for each $v\in[0,1]$, $t\mapsto X_t(v)$ is c\`{a}dl\`{a}g on
$[0,\infty)$ and solves (\ref{4.3});

\item
for each $t\ge0$, $v\mapsto X_t(v)$ is nondecreasing and
c\`{a}dl\`{a}g on $[0,1]$ with $X_t(0)\ge0$ and $X_t(1)\le1$.
\end{longlist}
\end{theorem}
\begin{pf}
This follows from Lemma \ref{l3.5} and Proposition \ref{p4.2}
by arguments as in the proof of Theorem \ref{t3.6}.
\end{pf}

We call the solution flow $\{X_t(v)\dvtx t\ge0, v\in[0,1]\}$ of (\ref{4.3})
specified in Theorem~\ref{t4.4} a \textit{generalized Fleming--Viot flow}
following Bertoin and Le Gall (\citeyear{BeL03}, \citeyear{BeL05},
\citeyear{BeL06}). The law of the flow is
determined by the parameters $(\sigma,b,\gamma,\nu)$.

Let $F[0,1]$ be the set of nondecreasing c\`{a}dl\`{a}g functions $f$ on
$[0,1]$ such that $0\le f(0)\le f(1)\le1$. Given a finite stopping time
$\tau$ and a function $\mu\in F[0,1]$, let $\{X^\mu_{\tau,t}(v)\dvtx
t\ge0\}$
be the solution of
%
%
\begin{eqnarray}\label{4.9}\quad
X^\mu_{\tau,t}(v) &=&\mu(v) + \int_\tau^{\tau+t}\int_0^1
\sigma\bigl[1_{\{u\le X^\mu_{\tau,s-}(v)\}} - X^\mu_{\tau,s-}(v)\bigr]
B(ds,du) \nonumber\\
&&{}
+ \int_\tau^{\tau+t} b[\gamma(v)-X^\mu_{\tau,s-}(v)] \,ds \\
&&{}
+ \int_\tau^{\tau+t}\int_0^1\int_0^1z\bigl[1_{\{u\le X^\mu_{\tau
,s-}(v)\}} -
X^\mu_{\tau,s-}(v)\bigr] \tilde{M}(ds,dz,du)\nonumber
\end{eqnarray}
and write simply $\{X^\mu_t(v)\dvtx t\ge0\}$ instead of $\{X^\mu
_{0,t}(v)\dvtx
t\ge0\}$. The pathwise uniqueness for the above equation follows from
that of (\ref{4.3}).\vadjust{\goodbreak} Let $F_{\tau,t}$ be the random operator on $F[0,1]$
that maps $\mu$ to $X^\mu_{\tau,t}$. As for the flow of
CBI-processes we
have
\begin{theorem}\label{t4.5} For any finite stopping time $\tau$ we have
$\mbf{P}\{X^\mu_{\tau+t} = F_{\tau,t}X^\mu_t$ for all $t\ge0\} = 1$.
\end{theorem}

For any sub-probability measure $\mu(dv)$ on $[0,1]$ with distribution
function $v\mapsto\mu(v)$, we write $X^\mu_t(dv)$ for the
random
sub-probability measure on $[0,1]$ determined by the random function
$v\mapsto X^\mu_t(v)$. We call $\{X_t^\mu\dvtx t\ge0\}$ the
\textit{generalized Fleming--Viot process} associated with the flow
$\{X_t^\mu(v)\dvtx t\ge0, v\in[0,1]\}$. The reader may refer to
\citet{Daw93} and \citet{EtK93} for the theory of classical
Fleming--Viot processes. To give some characterizations of the generalized
Fleming--Viot process, let us consider the step function
%
%
\begin{equation}\label{4.10}
f(u) = c_01_{\{0\}}(u) + \sum_{i=1}^n c_i1_{(a_{i-1},a_i]}(u),\qquad
u\in[0,1],
\end{equation}
where $\{c_0,c_1, \ldots, c_n\}\subset\mbb{R}$ and $\{0=a_0< a_1<
\cdots< a_n=1\}$ is a partition of $[0,1]$. For this function we
have
%
%
\begin{equation}\label{4.11}
\langle X^\mu_t,f\rangle
=
c_0X^\mu_t(0) + \sum_{i=1}^n c_i[X^\mu_t(a_i) - X^\mu_t(a_{i-1})].
\end{equation}
By (\ref{4.9}) and (\ref{4.11}) we have
%
%
\begin{eqnarray}\label{4.12}
\langle X^\mu_t,f\rangle
&=&
\langle\mu,f\rangle+ \int_0^t\int_0^1 \sigma[g^\mu_{s-}(u) -
\langle X^\mu_{s-},f\rangle] B(ds,du) \nonumber\\
&&{}
+ \int_0^t b[\langle\gamma,f\rangle- \langle X^\mu
_{s-},f\rangle] \,ds \\
&&{}
+ \int_0^t\int_0^1\int_0^1 z [g^\mu_{s-}(u) -
\langle X^\mu_{s-},f\rangle] \tilde{M}(ds,dz,du),\nonumber
\end{eqnarray}
where
%
%
\begin{equation}\label{4.13}
g^\mu_s(u) = c_01_{\{u\le X^\mu_s(0)\}} + \sum_{i=1}^n
c_i1_{\{X^\mu_s(a_{i-1})< u\le X^\mu_s(a_i)\}}.
\end{equation}
The proofs of the following three results are similar to those for
CBI-processes.

\begin{theorem}\label{t4.6} The generalized Fleming--Viot process
$\{X^\mu_t\dvtx t\ge0\}$ defined above is an almost surely c\`{a}dl\`{a}g
strong Markov process with $X^\mu_0 = \mu$.
\end{theorem}
\begin{prop}\label{p4.7} For any $t\ge0$ and $f\in B[0,1]$ we have
%
%
\begin{equation}\label{4.14}
\mbf{E}[\langle X^\mu_t,f\rangle]
=
\langle\mu,f\rangle e^{-bt} + \langle\gamma,f\rangle(1 - e^{-bt}).
\end{equation}
\end{prop}
\begin{theorem}\label{t4.8}
$\!\!\!\!$For any $f\!\in\! B[0,1]$ the process
$\{\langle X^\mu_t,f\rangle\dvtx t\,{\ge}\,0\}$ has a~c\`{a}dl\`{a}g
modification. Moreover,
there is a locally bounded function $t\mapsto C(t)$ so that
%
%
\begin{equation}\label{4.15}\qquad
\mbf{E}\Bigl[\sup_{0\le s\le t}\langle X^\mu_s,f\rangle\Bigr]
\le
C(t)[\langle\mu,f\rangle+ \langle\gamma,f\rangle
+ \langle\mu,f^2\rangle^{1/2} + \langle\gamma,f^2\rangle
^{1/2}]
\end{equation}
for any $t\ge0$ and $f\in B[0,1]^+$.
\end{theorem}

The generalized Fleming--Viot process can be characterized in terms of a~martingale problem. Given any finite family $\{f_1,\ldots,f_p\}\subset
B[0,1]$, write
%
%
\begin{equation}\label{4.16}
G_{p,\{f_i\}}(\eta) = \prod_{i=1}^p\langle\eta,f_i\rangle,\qquad
\eta\in M_1[0,1].
\end{equation}
Let $\mcr{D}_1(L)$ be the linear span of the functions on $M_1[0,1]$ of
the form (\ref{4.16}), and let $L$ be the linear operator on $\mcr{D}_1(L)$
defined by
%
%
\begin{eqnarray}\label{4.17}\quad\qquad
LG_{p,\{f_i\}}(\eta)
&=&
\sigma^2\sum_{i<j}\Biggl[\langle\eta,f_if_j\rangle\prod_{k\neq
i,j} \langle\eta,f_k\rangle-
\prod_{k=1}^p\langle\eta,f_k\rangle\Biggr] \nonumber\\[-2pt]
&&{}
+ \sum_{I\subset\{1,\ldots,p\},|I|\ge2} \beta_{p,|I|}\Biggl[
\biggl\langle\eta,
\prod_{i\in I} f_i\biggr\rangle\prod_{j\notin I}\langle\eta
,f_j\rangle- \prod_{k=1}^p
\langle\eta,f_k\rangle\Biggr] \\[-2pt]
&&{}
+ b\sum_{i=1}^p\Biggl[\langle\gamma,f_i\rangle\prod_{k\neq
i}\langle\eta,f_k\rangle-
\prod_{k=1}^p\langle\eta,f_k\rangle\Biggr],\nonumber
\end{eqnarray}
where $|I|$ denotes the cardinality of $I\subset\{1,\ldots,p\}$ and
\[
\beta_{p,|I|} = \int_0^1z^{|I|}(1-z)^{p-|I|}\nu(dz).
\]

\begin{theorem}\label{t4.9} The generalized Fleming--Viot process
$\{X^\mu_t\dvtx t\ge0\}$ is the unique solution of the following
martingale problem: for any $p\ge1$ and $\{f_1,\ldots, f_p\}\subset
B[0,1]$,
%
%
\begin{equation}\label{4.18}
G_{p,\{f_i\}}(X^\mu_t)
=
G_{p,\{f_i\}}(\mu) + \int_0^t LG_{p,\{f_i\}}(X^\mu_s)\,ds +
\mbox{mart.}
\end{equation}
\end{theorem}

\begin{pf}
We first consider a collection of step functions
$\{f_1,\ldots,f_p\}$. Let $g^\mu_i(s,u)$ be defined by (\ref{4.13}) with
$f=f_i$. Since the compensation of the Poisson random measure in
(\ref{4.12}) can be disregarded, by It\^{o}'s formula we
get\looseness=-1
\begin{eqnarray*}
&&G_{p,\{f_i\}}(X^\mu_t)\\[-2pt]
&&\qquad=
G_{p,\{f_i\}}(\mu) + \sigma^2\int_0^t ds\int_0^1 \biggl[\sum_{i<j}
h^\mu_i(s,u)h^\mu_j(s,u)
\prod_{k\neq i,j}\langle X^\mu_s,f_k\rangle\biggr]\, du \\[-2pt]
&&\qquad\quad{}
+ \int_0^tds\int_0^1 \nu(dz)\int_0^1 \Biggl\{\prod_{k=1}^p
[\langle X^\mu_s,f_k\rangle+ zh^\mu_k(s,u)]
- \prod_{k=1}^p \langle X^\mu_s,f_k\rangle\Biggr\}\, du \\[-2pt]
&&\qquad\quad{}
+ b\int_0^t \sum_{i=1}^p [\langle\gamma,f_i\rangle- \langle X
^\mu_s,f_i\rangle]
\prod_{k\neq i}\langle X^\mu_s,f_k\rangle\,ds + \mbox{mart.} \\[-2pt]
&&\qquad=
G_{p,\{f_i\}}(\mu) + \sigma^2\int_0^t ds\int_0^1 \biggl[\sum_{i<j}
l^\mu_i(u)l^\mu_j(u) \prod_{k\neq i,j}\langle X^\mu_s,f_k\rangle\biggr] X^\mu
_s(du) \\[-2pt]
&&\qquad\quad{}
+ \int_0^tds\int_0^1 \nu(dz)\int_0^1 \Biggl\{\prod_{k=1}^p
[\langle X^\mu_s,f_k\rangle+ zl^\mu_k(u)]
- \prod_{k=1}^p \langle X^\mu_s,f_k\rangle\Biggr\} X^\mu_s(du) \\[-2pt]
&&\qquad\quad{}
+ b\int_0^t \sum_{i=1}^p \Biggl[\langle\gamma,f_i\rangle\prod
_{k\neq i}
\langle X^\mu_s,f_k\rangle- \prod_{k=1}^p\langle X^\mu_s,f_k\rangle
\Biggr] \,ds + \mbox{mart.},\vspace*{-2pt}
\end{eqnarray*}\looseness=0
where $h^\mu_i(s,u) = g^\mu_i(s,u) - \langle X^\mu_s,f_i\rangle$ and
$l^\mu_i(u) = f_i(u) - \langle X^\mu_s,f_i\rangle$. It is simple to show
\[
\int_0^1 l^\mu_i(u)l^\mu_j(u) X^\mu_s(du)
=
\langle X^\mu_s,f_if_j\rangle- \langle X^\mu_s,f_i\rangle\langle X
^\mu_s,f_j\rangle.\vspace*{-2pt}
\]
Then we continue with
\begin{eqnarray*}
\hspace*{-4pt}&&G_{p,\{f_i\}}(X^\mu_t)\\[-2.5pt]
\hspace*{-4pt}&&\qquad=
G_{p,\{f_i\}}(\mu) + \sigma^2\int_0^t \sum_{i<j}
\Biggl[\langle X^\mu_s,f_if_j\rangle\prod_{k\neq i,j}\langle X^\mu
_s,f_k\rangle
- \prod_{k=1}^p\langle X^\mu_s,f_k\rangle\Biggr] \,ds \\[-2.5pt]
\hspace*{-4pt}&&\qquad\quad{}
+ \int_0^tds\int_0^1 \nu(dz)\int_0^1 \Biggl\{\prod_{k=1}^p
[(1-z)\langle X^\mu_s,f_k\rangle+ zf_k(u)]\\[-2.5pt]
\hspace*{-4pt}&&\hspace*{205.6pt}{} - \prod_{k=1}^p \langle X^\mu_s,f_k\rangle\Biggr\} X^\mu_s(du) \\[-2.5pt]
\hspace*{-4pt}&&\qquad\quad{}
+ b\int_0^t \sum_{i=1}^p \Biggl[\langle\gamma,f_i\rangle\prod
_{k\neq i}
\langle X^\mu_s,f_k\rangle- \prod_{k=1}^p\langle X^\mu_s,f_k\rangle
\Biggr] \,ds +
\mbox{mart.} \\[-2.5pt]
\hspace*{-4pt}&&\qquad=
G_{p,\{f_i\}}(\mu) + \sigma^2\int_0^t \sum_{i<j}
\Biggl[\langle X^\mu_s,f_if_j\rangle\prod_{k\neq i,j}\langle X^\mu
_s,f_k\rangle - \prod_{k=1}^p\langle X^\mu_s,f_k\rangle\Biggr] \,ds \\[-2.5pt]
\hspace*{-4pt}&&\qquad\quad{}
+ \int_0^tds\int_0^1 \nu(dz)\int_0^1 \Biggl\{\sum_{I\subset
\{1,\ldots,p\}}z^{|I|}(1-z)^{p-|I|}\prod_{i\in I} f_i(u)
\prod_{j\notin I}\langle X^\mu_s,f_j\rangle\\[-2.5pt]
\hspace*{-4pt}&&\hspace*{232pt}{}- \prod_{k=1}^p
\langle X^\mu_s,f_k\rangle\Biggr\}
X^\mu_s(du) \\[-2.5pt]
\hspace*{-4pt}&&\qquad\quad{}
+ b\int_0^t \sum_{i=1}^p \Biggl[\langle\gamma,f_i\rangle\prod
_{k\neq i}
\langle X^\mu_s,f_k\rangle- \prod_{k=1}^p\langle X^\mu_s,f_k\rangle
\Biggr] \,ds +
\mbox{mart.} \\[-2.5pt]
\hspace*{-4pt}&&\qquad=
G_{p,\{f_i\}}(\mu) + \sigma^2\int_0^t \sum_{i<j}
\Biggl[\langle X^\mu_s,f_if_j\rangle\prod_{k\neq i,j}\langle X^\mu
_s,f_k\rangle
- \prod_{k=1}^p\langle X^\mu_s,f_k\rangle\Biggr] \,ds \\[-2.5pt]
\hspace*{-4pt}&&\qquad\quad{}
+ \int_0^tds\int_0^1 \nu(dz)\int_0^1 \Biggl\{\sum_{I\subset
\{1,\ldots,p\}}z^{|I|}(1-z)^{p-|I|}\Biggl[\prod_{i\in I} f_i(u)
\prod_{j\notin I}\langle X^\mu_s,f_j\rangle\\[-2.5pt]
\hspace*{-4pt}&&\qquad\quad\hspace*{197pt}{} - \prod_{k=1}^p
\langle X^\mu_s,f_k\rangle\Biggr]\Biggr\}
X^\mu_s(du) \\[-2.5pt]
\hspace*{-4pt}&&\qquad\quad{}
+ b\int_0^t \sum_{i=1}^p \Biggl[\langle\gamma,f_i\rangle\prod
_{k\neq i}
\langle X^\mu_s,f_k\rangle- \prod_{k=1}^p\langle X^\mu_s,f_k\rangle
\Biggr] \,ds +
\mbox{mart.}
\end{eqnarray*}
That gives (\ref{4.18}) for step functions $\{f_1,\ldots,f_p\}$. For
$\{f_1,\ldots,f_p\}\subset B[0,1]$ one can show (\ref{4.18}) by
approximating the functions in the space $L^2(\mu+\gamma)$ using bounded
sequences of step functions. Since $\{X^\mu_t\dvtx t\ge0\}$ is a Markov
process, and $\mcr{D}_1(L)$ separates probability measures on $M[0,1]$,
the uniqueness for the martingale problem holds [see Ethier and Kurtz
(\citeyear{EtK86}), page~182].
\end{pf}

In particular, if $\mu(1) = \gamma(1) = 1$, we have $X_t^\mu(1) = 1$
for all $t\ge0$, and the corresponding generalized Fleming--Viot
process $\{X_t^\mu\dvtx t\ge0\}$ is a probability-valued Markov process
with generator $L$ defined by
%
%
\begin{eqnarray}\label{4.19}\qquad
LG_{p,\{f_i\}}(\eta)
&=&
\sigma^2\sum_{i<j}\Biggl[\langle\eta,f_if_j\rangle\prod_{k\neq
i,j} \langle\eta,f_k\rangle-
\prod_{k=1}^p\langle\eta,f_k\rangle\Biggr] \nonumber\\[1.2pt]
&&{}
+ \sum_{I\subset\{1,\ldots,p\},|I|\ge2} \beta_{p,|I|}\Biggl[
\biggl\langle\eta,
\prod_{i\in I} f_i\biggr\rangle\prod_{j\notin I}\langle\eta
,f_j\rangle- \prod_{k=1}^p
\langle\eta,f_k\rangle\Biggr] \\[1.2pt]
&&{}
+ \sum_{i=1}^p\langle\eta,Af_i\rangle\prod_{k\neq i}\langle\eta
,f_k\rangle,\nonumber
\end{eqnarray}
where
\[
Af(x) = b\int_{[0,1]}[f(y)-f(x)]\gamma(dy),\qquad x\in[0,1].
\]
This is a generalization of a classical Fleming--Viot process [see, e.g.,
Ethier and Kurtz (\citeyear{EtK93}), page 351]. On the other hand, for
$b=0$ the solution
flow $\{X_t^\mu(v)\dvtx t\ge0, 0\le v\le1\}$ of (\ref{4.3})
corresponds to\vspace*{1pt}
the $\Lambda$-coalescent process with $\Lambda(dz) = \sigma
^2\delta_0
+ z^2\nu(dz)$, which is clear from (\ref{4.18}) and the\vadjust{\goodbreak} martingale problem
given by Theorem 1 in \citet{BeL05}. For $b>0$ it seems the
flow determines a coalescent process with a spatial structure. A~serious
exploration in the subject would be of interest to the understanding of
the related dynamic systems.


\section{Scaling limit theorems}\label{sec5}

In this section, we prove some limit theorems for the generalized
Fleming--Viot flows. We shall present the results in the setting of
measure-valued processes and through the use of Markov process
arguments. These are
different from the approach of \citet{BeL06}, who used the
analysis of characteristics of semimartingales. For each $k\ge1$ let
$\sigma_k\ge0$ and $b_k\ge0$ be two constants, let $z^2\nu_k(dz)$
be a
finite measure on $(0,1]$ and let $v\mapsto\gamma_k(v)$ be a
nondecreasing continuous function on $[0,1]$ so that $0\le
\gamma_k(v)\le1$ for all $0\le v\le1$. We denote by $\gamma_k(dv)$ the
sub-probability measure on $[0,1]$ so that $\gamma_k([0,v]) =
\gamma_k(v)$ for $0\le v\le1$. Let $\{X_t^k(v)\dvtx t\ge0,v\in[0,1]\}$ be
a generalized Fleming--Viot flow with parameters
$(\sigma_k,b_k,\gamma_k,\nu_k)$ and with $X_0^k(v) = v$ for $v\in[0,1]$.
Let $Y_k(t,v) = kX_{kt}^k(k^{-1}v)$ for $t\ge0$ and $v\in[0,k]$. Let
$\eta_k(z) = k\gamma_k(k^{-1}z)$ and $m_k(dz) = \nu_k(k^{-1}dz)$ for
$z\in(0,k]$. In view of (\ref{4.3}), we can also define $\{Y_k(t,v)\dvtx
t\ge0,v\in[0,k]\}$ directly by
%
%
\begin{eqnarray}\label{5.1}\quad
Y_k(t,v)
&=&
v + k\sigma_k\int_0^t\int_0^k \bigl[1_{\{u\le Y_k(s-,v)\}} -
k^{-1}Y_k(s-,v)\bigr] W_k(ds,du) \nonumber\\[-2pt]
&&{}
+ kb_k\int_0^t [\eta_k(v) - Y_k(s-,v)] \,ds \\[-2pt]
&&{}
+ \int_0^t\int_0^k\int_0^k z\bigl[1_{\{u\le Y_k(s-,v)\}} -
k^{-1}Y_k(s-,v)\bigr] \tilde{N}_k(ds,dz,du),\nonumber
\end{eqnarray}
where $\{W_k(ds,du)\}$ is a white noise on $(0,\infty)\times(0,k]$ with
intensity $ds\,du$, and $\{N_k(ds,dz,du)\}$ is a Poisson random measure on
$(0,\infty)\times(0,k]^2$ with intensity $ds\,m_k(dz)\,du$. In the
sequel, we
assume $k\ge a$ for fixed a constant $a\ge0$. Then the rescaled flow
$\{Y_k(t,v)\dvtx t\ge0,v\in[0,k]\}$ induces an $M[0,a]$-valued process
$\{Y_k^a(t)\dvtx t\ge0\}$. We are interested in the asymptotic
behavior of
$\{Y_k^a(t)\dvtx t\ge0\}$ as $k\to\infty$. Recall that $\lambda$
denotes the
Lebesgue measure on $[0,\infty)$.
\begin{lemma}\label{l5.1} For any $G\in C^2(\mbb{R})$ and $f\in
C[0,a]$ we
have
\begin{eqnarray*}
&&G(\langle Y_k^a(t),f\rangle)\\[-2pt]
&&\qquad=
G(\langle\lambda,f\rangle) + kb_k\int_0^t G^\prime(\langle Y
_k^a(s),f\rangle)
\langle\eta_k,f\rangle\,ds \\[-2pt]
&&\qquad\quad{}
- kb_k\int_0^t G^\prime(\langle Y_k^a(s),f\rangle)\langle Y
_k^a(s),f\rangle\,ds \\[-2pt]
&&\qquad\quad{}
+ \frac{1}{2}k^2\sigma_k^2\int_0^t G^{\prime\prime}(\langle Y
_k^a(s),f\rangle)
\langle Y_k^a(s),f^2\rangle\,ds \\[-2pt]
&&\qquad\quad{}
- \frac{1}{2}k\sigma_k^2\int_0^t G^{\prime\prime}(\langle Y
_k^a(s),f\rangle)
\langle Y_k^a(s),f\rangle^2 \,ds \\[-2pt]
&&\qquad\quad{}
+ \int_0^tds\int_0^km_k(dz)\int_{[0,a]}\bigl\{G\bigl(\langle Y
_k^a(s),f\rangle+
zf(x)\bigr)
- G(\langle Y_k^a(s),f\rangle)\\[-2pt]
&&\qquad\quad\hspace*{145pt}{} - G^\prime(\langle Y_k^a(s),f\rangle
)zf(x)\bigr\}
Y_k^a(s,dx) \\[-2pt]
&&\qquad\quad{}
+ \int_0^tds\int_0^k [\varepsilon_k(s,z)+\xi_k(s,z)]m_k(dz) +
\mbox{local mart.},
\end{eqnarray*}
where
\begin{eqnarray*}
\varepsilon_k(s,z)
&=&
\int_0^k \bigl\{G\bigl(\langle Y_k^a(s),f\rangle+ z[f(x) - k^{-1} \langle Y
_k^a(s),f\rangle]\bigr)
\\[-2pt]
&&\hspace*{17.5pt}{}
- G\bigl(\langle Y_k^a(s),f\rangle+ zf(x)\bigr) \\[-2pt]
&&\hspace*{17.5pt}{}
- k^{-1}G^\prime(\langle Y_k^a(s),f\rangle)z\langle Y
_k^a(s),f\rangle\bigr\} Y_k^a(s,dx)
\end{eqnarray*}
and
\begin{eqnarray*}
\xi_k(s,z)
&=&
[k-Y_k(s,a)]\\[-2pt]
&&{}\times\bigl[G\bigl(\langle Y_k^a(s),f\rangle- k^{-1}z\langle Y
_k^a(s),f\rangle\bigr) \\[-2pt]
&&\hspace*{13.8pt}{}
- G(\langle Y_k^a(s),f\rangle) + k^{-1}G^\prime(\langle Y
_k^a(s),f\rangle)
z\langle Y_k^a(s),f\rangle\bigr].
\end{eqnarray*}
\end{lemma}
\begin{pf}
For the step function defined by (\ref{3.13}) we get from
(\ref{5.1}) that
%
%
\begin{eqnarray}\label{5.2}
\langle Y_k^a(t),f\rangle
&=&
\langle\lambda,f\rangle+ k\sigma_k\int_0^t\int_0^k h_k(s-,u)
W_k(ds,du) \nonumber\\[-2pt]
&&{}
+ kb_k\int_0^t [\langle\eta_k,f\rangle- \langle Y
_k^a(s-),f\rangle] \,ds \\[-2pt]
&&{}
+ \int_0^t\int_0^k\int_0^k zh_k(s-,u) \tilde{N}_k(ds,dz,du),\nonumber
\end{eqnarray}
where $h_k(s,u) = g_k(s,u) - k^{-1}\langle Y_k^a(s),f\rangle$ and
%
%
\begin{equation}\label{5.3}
g_k(s,u) = c_0 1_{\{u\le Y_k(s,0)\}} + \sum_{i=1}^n c_i
1_{\{Y_k(s,a_{i-1})< u\le Y_k(s,a_i)\}}.
\end{equation}
Let $l_k(s,x) = f(x) - k^{-1}\langle Y_k^a(s),f\rangle$. By (\ref
{5.2}) and
It\^{o}'s formula,
\begin{eqnarray*}
&&G(\langle Y_k^a(t),f\rangle)\\
&&\qquad=
G(\langle\lambda,f\rangle) + kb_k\int_0^t G^\prime(\langle Y
_k^a(s),f\rangle)
[\langle\eta_k,f\rangle- \langle Y_k^a(s),f\rangle] \,ds \\
&&\qquad\quad{}
+ \frac{1}{2}k^2\sigma_k^2\int_0^t G^{\prime\prime}
(\langle Y_k^a(s),f\rangle) \,ds\int_0^k h_k(s,u)^2 \,du \\
&&\qquad\quad{}
+ \int_0^tds\int_0^k m_k(dz)\int_0^k \bigl\{G\bigl(\langle Y
_k^a(s),f\rangle+
zh_k(s,u)\bigr) \\
&&\qquad\quad\hspace*{108.7pt}{}
- G(\langle Y_k^a(s),f\rangle) \\
&&\qquad\quad\hspace*{108.7pt}{}- G^\prime(\langle Y_k^a(s),f\rangle
)zh_k(s,u)\bigr\}\, du
\nonumber\\
&&\qquad\quad{}
+ \mbox{local mart.} \nonumber\\
&&\qquad=
G(\langle\lambda,f\rangle) + kb_k\int_0^t G^\prime(\langle Y
_k^a(s),f\rangle)
[\langle\eta_k,f\rangle- \langle Y_k^a(s),f\rangle] \,ds \\
&&\qquad\quad{}
+ \frac{1}{2}k^2\sigma_k^2\int_0^t G^{\prime\prime}(\langle Y
_k^a(s),f\rangle)
[\langle Y_k^a(s),f^2\rangle- k^{-1}\langle Y_k^a(s),f\rangle] \,ds \\
&&\qquad\quad{}
+ \int_0^tds\int_0^km_k(dz)\int_{[0,a]}\bigl\{G\bigl(\langle Y
_k^a(s),f\rangle+
zl_k(s,x)\bigr) \\
&&\qquad\quad\hspace*{116.5pt}{}
- G(\langle Y_k^a(s),f\rangle) \\
&&\qquad\quad\hspace*{116.5pt}{}- G^\prime(\langle Y_k^a(s),f\rangle
)zl_k(s,x)\bigr\}
Y_k^a(s,dx) \\
&&\qquad\quad{}
+ \int_0^t[k-Y_k(s,a)]\,ds\int_0^k \bigl\{G\bigl(\langle Y_k^a(s),f\rangle-
k^{-1}z\langle Y_k^a(s),f\rangle\bigr) \\
&&\qquad\quad\hspace*{118.5pt}{}
- G(\langle Y_k^a(s),f\rangle)\\
&&\qquad\quad\hspace*{118.5pt}{} +
k^{-1}G^\prime(\langle Y_k^a(s),f\rangle)z\langle Y_k^a(s),f\rangle
\bigr\} m_k(dz) \nonumber\\
&&\qquad\quad{}
+ \mbox{local mart.}
\end{eqnarray*}
That gives the desired result for the step function. For $f\in C[0,a]$ it
follows by approximating the function by a sequence of step functions.
\end{pf}
\begin{lemma}\label{l5.2} For $t\ge0$ and $f\in C[0,a]^+$ we have
\begin{eqnarray*}
&&\mbf{E}\Bigl[\sup_{0\le s\le t}\langle Y_k^a(s),f\rangle
\Bigr] \\
&&\qquad
\le\langle\lambda,f\rangle+ kb_k\langle\eta_k,f\rangle t+ 4t
[\langle\lambda,f\rangle+
\langle\eta_k,f\rangle]\int_1^k z m_k(dz) \\
&&\qquad\quad{}
+ 2\sqrt{t}[\langle\lambda,f^2\rangle+\langle\eta
_k,f^2\rangle]^{1/2}
\biggl[\sigma+ \biggl(\int_0^1 z^2m_k(dz)\biggr)^{1/2}\biggr].
\end{eqnarray*}
\end{lemma}
\begin{pf}
We first consider a nonnegative step function given by
(\ref{3.13}) with $\{c_0,c_1,\ldots,c_n\}\subset\mbb{R}_+$.
Let
$g_k(s,u)$ and $h_k(s,u)$ be defined as in the proof of Lemma \ref{l5.1}.
By (\ref{5.2}) and Doob's martingale inequality we get
\begin{eqnarray*}
&&\mbf{E}\Bigl[\sup_{0\le s\le t}\langle Y_k^a(s),f\rangle
\Bigr] \\[-2pt]
&&\qquad
\le\langle\lambda,f\rangle+ 2k\sigma_k\mbf{E}^{1/2}
\biggl\{\biggl[\int_0^t\int_0^k h_k(s-,u) W(ds,du)\biggr]^2\biggr\}
\\[-2pt]
&&\qquad\quad{}
+ kb_k\langle\eta_k,f\rangle t+ \mbf{E}\biggl[\int_0^tds\int_1^k
z m_k(dz) \int_0^k
|h_k(s-,u)| \,du\biggr] \\[-2pt]
&&\qquad\quad{}
+ \mbf{E}\biggl[\int_0^t\int_1^k\int_0^k z|h_k(s-,u)|
N_k(ds,dz,du)\biggr] \\[-2pt]
&&\qquad\quad{}
+ 2\mbf{E}^{1/2}\biggl\{\biggl[\int_0^t\int_0^1\int_0^k z
h_k(s-,u) \tilde{N}_k(ds,dz,du)\biggr]^2\biggr\}.
\end{eqnarray*}
It then follows that
\begin{eqnarray*}
&&\mbf{E}\Bigl[\sup_{0\le s\le t}\langle Y_k^a(s),f\rangle
\Bigr] \\[-2pt]
&&\qquad
\le\langle\lambda,f\rangle+ 2k\sigma_k \mbf{E}^{1/2}
\biggl\{\int_0^tds
\int_0^k h_k(s,u)^2 \,du\biggr\} \\[-2pt]
&&\qquad\quad{}
+ kb_k\langle\eta_k,f\rangle t+ 2\mbf{E}\biggl\{\int_0^tds\int
_1^k z
m_k(dz)\int_0^k |h_k(s,u)| \,du\biggr\} \\[-2pt]
&&\qquad\quad{}
+ 2\mbf{E}^{1/2}\biggl\{\int_0^tds \int_0^1 z^2 m_k(dz)
\int_0^k h_k(s,u)^2 \,du\biggr\} \\[-2pt]
&&\qquad
\le\langle\lambda,f\rangle+ kb_k\langle\eta_k,f\rangle t+ 4\mbf
{E}\biggl[\int_0^t
\langle Y_k^a(s),f\rangle\,ds\int_1^k z m_k(dz)\biggr] \\[-2pt]
&&\qquad\quad{}
+ 2\mbf{E}^{1/2}\biggl[\int_0^t \langle Y
_k^a(s),f^2\rangle\,
ds\biggr]\biggl[k\sigma_k + \biggl(\int_0^1 z^2
m_k(dz)\biggr)^{1/2}\biggr].
\end{eqnarray*}
By Proposition \ref{p4.7} one can see
\[
\mbf{E}[\langle Y_k^a(t),f\rangle]
=
\langle\lambda,f\rangle e^{-kb_kt} + \langle\eta_k,f\rangle(1-e^{-kb_kt})
\le
\langle\lambda,f\rangle+ \langle\eta_k,f\rangle.
\]
Then we have the desired inequality for the step function. The
inequality for $f\in C[0,a]^+$ follows by approximating this
function with a bounded sequence of positive step functions.
\end{pf}
\begin{lemma}\label{l5.3} Let $\tau_k$ be a bounded stopping time for
$\{Y_k^a(t)\dvtx t\ge0\}$. Then for any $t\ge0$ and $f\in C[0,a]$ we
have
%
%
\begin{eqnarray}\label{5.4}
&&\mbf{E}\{|\langle Y_k^a(\tau_k+t),f\rangle-
\langle Y_k^a(\tau_k),f\rangle|\} \nonumber\\
&&\qquad
\le\mbf{E}^{1/2}\biggl[\int_0^t \langle Y_k^a(\tau_k+s),
f^2\rangle\,
ds\biggr] \biggl[k\sigma_k + \biggl(\int_0^1 z^2
m_k(dz)\biggr)^{1/2}\biggr] \nonumber\\[-10pt]\\[-12pt]
&&\qquad\quad{}
+ kb_k\mbf{E}\biggl[\int_0^t \bigl(\langle\eta_k,|f|\rangle+
\langle Y_k^a(\tau_k+s),|f|\rangle\bigr) \,ds\biggr] \nonumber\\[-2.5pt]
&&\qquad\quad{}
+ 4\mbf{E}\biggl[\int_0^t \langle Y_k^a(\tau_k+s),|f|\rangle\,
ds\int_1^k z
m_k(dz)\biggr].\nonumber\vspace*{-2pt}
\end{eqnarray}
\end{lemma}

\begin{pf}
We first consider the step function given by (\ref{3.13}). Let
$g_k(s,u)$ and $h_k(s,u)$ be defined as in the proof of Lemma \ref{l5.1}.
{From} (\ref{5.2}) we have
\begin{eqnarray*}
&&\mbf{E}\{|\langle Y_k^a(\tau_k+t),f\rangle-
\langle Y_k^a(\tau_k),f\rangle|\} \\[-2.5pt]
&&\qquad
\le k\sigma_k\mbf{E}^{1/2}\biggl\{\biggl[\int_0^t\int_0^k
h_k(\tau_k+s-,u) W(\tau_k+ds,du)\biggr]^2\biggr\} \\[-2.5pt]
&&\qquad\quad{}
+ kb_k\mbf{E}\biggl[\int_0^t |\langle\eta_k,f\rangle-
\langle Y_k^a(\tau_k+s-),f\rangle| \,ds\biggr] \\[-2.5pt]
&&\qquad\quad{}
+ \mbf{E}^{1/2}\biggl\{\biggl[\int_0^t\int_0^1\int_0^k
zh_k(\tau_k+s-,u) \tilde{N}_k(\tau_k+ds,dz,du)\biggr]^2\biggr\} \\[-2.5pt]
&&\qquad\quad{}
+ \mbf{E}\biggl[\int_0^t\int_1^k\int_0^k z |h_k(\tau_k+s-,u)|
N_k(\tau_k+ds,dz,du)\biggr] \\[-2.5pt]
&&\qquad\quad{}
+ \mbf{E}\biggl[\int_0^tds\int_1^k z m_k(dz) \int_0^k
|h_k(\tau_k+s-,u)| \,du\biggr].\vspace*{-2pt}
\end{eqnarray*}
By the property of independent increments of the white noise and the
Poisson random measure,
\begin{eqnarray*}
&&\mbf{E}\{|\langle Y_k^a(\tau_k+t),f\rangle-
\langle Y_k^a(\tau_k),f\rangle|\} \\[-2.5pt]
&&\qquad
\le k\sigma_k\mbf{E}^{1/2}\biggl\{\int_0^tds\int_0^k
h_k(\tau_k+s,u)^2 \,du\biggr\} \\[-2.5pt]
&&\qquad\quad{}
+ kb_k\mbf{E}\biggl[\int_0^t \bigl(\langle\eta_k,|f|\rangle+
\langle Y_k^a(\tau_k+s),|f|\rangle\bigr) \,ds\biggr] \\[-2.5pt]
&&\qquad\quad{}
+ \mbf{E}^{1/2}\biggl\{\int_0^tds \int_0^1 z^2 m_k(dz)
\int_0^k
h_k(\tau_k+s,u)^2 \,du\biggr\} \\[-2.5pt]
&&\qquad\quad{}
+ 2\mbf{E}\biggl[\int_0^tds\int_1^k z m_k(dz) \int_0^k
|h_k(\tau_k+s,u)|\,
du\biggr] \\[-2.5pt]
&&\qquad
\le\mbf{E}^{1/2}\biggl[\int_0^t\langle Y_k^a(\tau
_k+s),f^2\rangle\,
ds\biggr]\biggl[k\sigma_k + \biggl(\int_0^1 z^2
m_k(dz)\biggr)^{1/2}\biggr] \\[-2.5pt]
&&\qquad\quad{}
+ kb_k\mbf{E}\biggl[\int_0^t \bigl(\langle\eta_k,|f|\rangle+
\langle Y_k^a(\tau_k+s),|f|\rangle\bigr) \,ds\biggr] \\[-2.5pt]
&&\qquad\quad{}
+ 4\mbf{E}\biggl[\int_0^t \langle Y_k^a(\tau_k+s),|f|\rangle
\,ds\int_1^k z
m_k(dz)\biggr].\vspace*{-2pt}\vadjust{\goodbreak}
\end{eqnarray*}
Then (\ref{5.4}) holds for the step function. For $f\in C[0,a]$ the
inequality follows by an approximation argument.
\end{pf}
\begin{lemma}\label{l5.4} Suppose that\vspace*{1pt} $kb_k\to b$, $\eta_k\to\eta$
weakly on $[0,a]$ and $k^2\sigma_k^2\times\delta_0(dz) + (z\land
z^2)m_k(dz)$ converges weakly on $[0,\infty)$ to a finite measure
$\sigma^2\times\delta_0(dz) + (z\land z^2)m(dz)$ as $k\to\infty$. Let
$\{0\le a_1< \cdots< a_n\}$ be an ordered set of constants. Then
$\{(Y_k^{a_1}(t), \ldots, Y_k^{a_n}(t))\dvtx t\ge0\}$, $k=1,2,\ldots$
is a tight sequence in $D([0,\infty),M[0,a_1]\times\cdots\times
M[0,a_n])$.
\end{lemma}
\begin{pf}
Let $\tau_k$ be a bounded stopping time for $\{Y_k^a(t)\dvtx
t\ge
0\}$ and assume the sequence $\{\tau_k\dvtx k=1,2,\ldots\}$ is uniformly
bounded. Let $f_i\in C[0,a_i]$ for $i=1,\ldots,n$. By (\ref{5.4}) we
see
%
%
\begin{eqnarray}\label{5.5}
&&\mbf{E}\Biggl\{\sum_{i=1}^n|\langle Y_k^{a_i}(\tau
_k+t),f_i\rangle-
\langle Y_k^{a_i}(\tau_k),f_i\rangle|\Biggr\} \nonumber\\
&&\qquad
\le\sum_{i=1}^n\mbf{E}^{1/2}\biggl[\int_0^t
\langle Y_k^{a_i}(\tau_k+s), f_i^2\rangle\,ds\biggr] \biggl[k\sigma_k +
\biggl(\int_0^1 z^2 m_k(dz)\biggr)^{1/2}\biggr] \nonumber\\[-8pt]\\[-8pt]
&&\qquad\quad{}
+ kb_k\sum_{i=1}^n\mbf{E}\biggl[\int_0^t \bigl(\langle\eta
_k,|f_i|\rangle+
\langle Y_k^{a_i}(\tau_k+s),|f_i|\rangle\bigr) \,ds\biggr] \nonumber\\
&&\qquad\quad{}
+ 4\sum_{i=1}^n\mbf{E}\biggl[\int_0^t \langle Y_k^{a_i}(\tau
_k+s),|f_i|\rangle\,
ds\int_1^k z m_k(dz)\biggr].\nonumber
\end{eqnarray}
Then the inequality in Lemma \ref{l5.2} implies
\[
\lim_{t\to0}\sup_{k\ge1} \mbf{E}\Biggl\{\sum_{i=1}^n|
\langle Y_k^{a_i}(\tau_k+t),f_i\rangle- \langle Y_k^{a_i}(\tau
_k),f_i\rangle|\Biggr\} = 0.
\]
By a criterion of \citet{Ald78}, the sequence
$\{(\langle Y_k^{a_1}(t),f_1\rangle, \ldots, \langle Y
_k^{a_n}(t),f_n\rangle)\dvtx\allowbreak t\ge0\}$
is tight in $D([0,\infty),\mbb{R}^n)$ [see also Ethier and Kurtz
(\citeyear{EtK86}), pages 137 and 138]. Then a simple extension of the
tightness
criterion of Roelly-Coppoletta (\citeyear{Roe86}) implies $\{(Y_k^{a_1}(t),
\ldots,
Y_k^{a_n}(t))\dvtx t\ge0\}$ is tight in $D([0,\infty),\allowbreak M[0,a_1]\times
\cdots\times M[0,a_n])$.
\end{pf}

Suppose that $\sigma\ge0$ and $b\ge0$ are two constants, $v\mapsto
\eta(v)$ is a nonnegative and nondecreasing continuous function on
$[0,\infty)$ and $(z\land z^2)m(dz)$ is a~finite measure on
$(0,\infty)$. Let $\eta(dv)$ be the Radon measure on $[0,\infty)$ so
that $\eta([0,v]) = \eta(v)$ for $v\ge0$. Suppose that
$\{W(ds,du)\}$ is a white noise on $(0,\infty)^2$ with intensity
$ds\,dz$ and $\{N(ds,dz,du)\}$ is a Poisson random measure on
$(0,\infty)^3$ with intensity $ds\,m(dz)\,du$. Let $\{X_t(v)\dvtx t\ge0,
v\ge0\}$ be the solution flow of the stochastic equation
%
%
\begin{eqnarray}\label{5.6}
X_t(v) &=& v + \sigma\int_0^t\int_0^{X_{s-}(v)} W(ds,du) +
b\int_0^t [\eta(v)-X_{s-}(v)] \,ds \nonumber\\[-8pt]\\[-8pt]
&&{}
+ \int_0^t \int_0^\infty\int_0^{X_{s-}(v)}
z\tilde{N}(ds,dz,du).\nonumber
\end{eqnarray}
By Theorem \ref{t3.11}, for each $a\ge0$ the flow $\{X_t(v)\dvtx t\ge
0, v\ge0\}$ induces an $M[0,a]$-valued immigration superprocess
$\{X^a_t\dvtx t\ge0\}$ which is the unique solution of the following
martingale problem: for every $G\in C^2(\mbb{R})$ and $f\in C[0,a]$,
%
%
\begin{eqnarray}\label{5.7}
&&G(\langle X_t,f\rangle)\nonumber\\
&&\qquad=
G(\langle\lambda,f\rangle) + b\int_0^t
G^\prime(\langle X_s,f\rangle)[\langle\eta,f\rangle-\langle X
_s,f\rangle] \,ds \nonumber\\
&&\qquad\quad
{}+ \frac{1}{2}\sigma^2 \int_0^t G^{\prime\prime}(\langle X
_s,f\rangle)
\langle X_s,f^2\rangle\,ds\nonumber\\[-8pt]\\[-8pt]
&&\qquad\quad{}
+ \int_0^tds\int_0^\infty m(dz)\int_{[0,a]} \bigl[G\bigl(\langle X
_s,f\rangle+
zf(x)\bigr) \nonumber\\
&&\qquad\quad\hspace*{115.1pt}{}
- G(\langle X_s,f\rangle) - G^\prime(\langle X_s,f\rangle
)zf(x)\bigr] X_s(dx) \nonumber\\
&&\qquad\quad{}
+ \mbox{local mart.}\nonumber
\end{eqnarray}
\begin{theorem}\label{t5.5} Suppose that $kb_k\to b$, $\eta_k\to\eta$
weakly on $[0,a]$ and\break $k^2\sigma_k^2\delta_0(dz) + (z\land
z^2)m_k(dz)$ converges weakly on $[0,\infty)$ to a finite measure
$\sigma^2\delta_0(dz) + (z\land z^2)m(dz)$ as $k\to\infty$. Then
$\{Y_k^a(t)\dvtx t\ge0\}$ converges to the immigration superprocess
$\{X_t^a\dvtx t\ge0\}$ in distribution on $D([0,\infty), M[0,a])$.
\end{theorem}

For the proof of the above theorem, let us make some
preparations. Since the solution of the martingale problem
(\ref{5.7}) is unique, it suffices to prove any weak limit point
$\{Z_t^a\dvtx t\ge0\}$ of the sequence $\{Y_k^a(t)\dvtx t\ge0\}$ is the
solution of the martingale problem. To simplify the notation we pass
to a~subsequence and simply assume $\{Y_k^a(t)\dvtx t\ge0\}$ converges
to $\{Z_t^a\dvtx t\ge0\}$ in distribution. Using Skorokhod's
representation theorem, we can also assume \mbox{$\{Y_k^a(t)\dvtx t\ge0\}$}
and $\{Z_t^a\dvtx t\ge0\}$ are defined on the same probability space
and $\{Y_k^a(t)\dvtx t\ge0\}$ converges a.s. to $\{Z_t^a\dvtx t\ge0\}
$ in
the topology of $D([0,\infty), M[0,a])$. For $n\ge1$ let
\[
\tau_n = \inf\biggl\{t\ge0\dvtx \sup_{k\ge1}\int_0^t [1+\langle Y_k^a(s)+Z^a_s,
1\rangle^2] \,ds\ge n\biggr\}.
\]
It is easy to see that $\tau_n\to\infty$ as $n\to\infty$.
\begin{lemma}\label{l5.6} Suppose that\vspace*{1pt} $kb_k\to b$, $\eta_k\to\eta$
weakly on $[0,a]$ and $k^2\sigma_k^2\times\delta_0(dz) + (z\land
z^2)m_k(dz)$ converges weakly on $[0,\infty)$ to a finite measure
$\sigma^2\times\delta_0(dz) + (z\land z^2)m(dz)$ as $k\to\infty$. Let
$\varepsilon_k(s,z)$ be defined as in Lemma \ref{l5.1}. Then for each
$n\ge1$ we have
\[
\mbf{E}\biggl[\int_0^{t\land\tau_n} ds\int_0^k |\varepsilon_k(s,z)|
m_k(dz)\biggr] \to0,\qquad k\to\infty.
\]
\end{lemma}
\begin{pf}
By the mean-value theorem, we have
\begin{eqnarray*}
\varepsilon_k(s,z)
&=&
\frac{1}{k} z\langle Y_k^a(s),f\rangle\\
&&{}\times\int_0^k\bigl[G^\prime
\bigl(\langle Y_k^a(s),f\rangle+
z\theta_k(s,x)\bigr)
- G^\prime(\langle Y_k^a(s),f\rangle)\bigr] Y_k^a(s,dx),
\end{eqnarray*}
where $\theta_k(s,x)$ takes values between $f(x)$ and $f(x) -
k^{-1}\langle Y_k^a(s),f\rangle$. Consequently,
\[
|\varepsilon_k(s,z)|
\le
\frac{2}{k}\|G^\prime\|z\langle Y_k^a(s),|f|\rangle\langle Y
_k^a(s),1\rangle
\le
\frac{2}{k}\|G^\prime\|\|f\|z\langle Y_k^a(s),1\rangle^2.
\]
Moreover, since $\langle Y_k^a(s),1\rangle\le k$, we get
\begin{eqnarray*}
|\varepsilon_k(s,z)|
&\le&
\frac{1}{k}\|G^{\prime\prime}\|z^2\langle Y_k^a(s),|f|\rangle\int_0^k
|\theta_k(s,x)| Y_k^a(s,dx) \\
&\le&
\frac{1}{k}\|G^{\prime\prime}\|z^2\langle Y_k^a(s),|f|\rangle\int
_0^k [|f(x)|
+ k^{-1}\langle Y_k^a(s),|f|\rangle] Y_k^a(s,dx) \\
&\le&
\frac{2}{k}\|f\|^2\|G^{\prime\prime}\| z^2\langle Y_k^a(s),1\rangle^2.
\end{eqnarray*}
It follows that
\begin{eqnarray*}
&&\mbf{E}\biggl[\int_0^{t\land\tau_n} ds\int_0^k |\varepsilon_k(s,z)|
m_k(dz)\biggr] \\
&&\qquad
\le\frac{C}{k}\int_0^k (z\land z^2) m_k(dz)
\mbf{E}\biggl[\int_0^{t\land\tau_n} \langle Y_k^a(s),1\rangle
^2\,ds\biggr] \\
&&\qquad
\le\frac{nC}{k}\int_0^k (z\land z^2) m_k(dz),
\end{eqnarray*}
where $C = 2\|f\|(\|G^\prime\|+\|G^{\prime\prime}\|\|f\|)$. The
right-hand side goes to zero as $k\to\infty$.
\end{pf}
\begin{lemma}\label{l5.7} Suppose\vspace*{1pt} that $kb_k\to b$, $\eta_k\to\eta$
weakly on $[0,a]$ and $k^2\sigma_k^2\times\delta_0(dz) + (z\land
z^2)m_k(dz)$ converges weakly on $[0,\infty)$ to a finite measure
$\sigma^2\times\delta_0(dz) + (z\land z^2)m(dz)$ as $k\to\infty$. Let
$\xi_k(s,z)$ be defined as in Lemma~\ref{l5.1}. Then for each $n\ge
1$ we have
\[
\mbf{E}\biggl[\int_0^{t\land\tau_n} ds\int_0^k |\xi_k(s,z)|
m_k(dz)\biggr] \to0,\qquad k\to\infty.
\]
\end{lemma}
\begin{pf}
It is elementary to see that
\begin{eqnarray*}
|\xi_k(s,z)|
&\le&
k\bigl|G\bigl(\langle Y_k^a(s),f\rangle- k^{-1}z\langle Y_k^a(s),f\rangle\bigr)
\\
&&\hspace*{6.6pt}{}- G(\langle Y_k^a(s),f\rangle)
+ k^{-1}G^\prime(\langle Y_k^a(s),f\rangle) z\langle Y
_k^a(s),f\rangle\bigr| \\
&\le&
\min\biggl\{2\|G^\prime\|z\langle Y_k^a(s),|f|\rangle,
\frac{1}{2k}\|G^{\prime\prime}\|z^2\langle Y_k^a(s),|f|\rangle
^2\biggr\} \nonumber\\
&\le&
C[1+\langle Y_k^a(s),1\rangle^2](z\land k^{-1}z^2),
\end{eqnarray*}
where $C = \|f\|(2\|G^\prime\|+\|f\|\|G^{\prime\prime}\|/2)$. Then we
have
\begin{eqnarray*}
&&\mbf{E}\biggl[\int_0^{t\land\tau_n} ds\int_0^k |\xi_k(s,z)|
m_k(dz)\biggr] \\
&&\qquad
\le C\int_0^k (z\land k^{-1}z^2) m_k(dz)
\mbf{E}\biggl\{\int_0^{t\land\tau_n} [1 + \langle Y_k^a(s),1\rangle
^2]\,ds\biggr\} \\
&&\qquad
\le nC\int_0^k (z\land k^{-1}z^2) m_k(dz).
\end{eqnarray*}
The right-hand side tends to zero as $k\to\infty$.
\end{pf}
\begin{pf*}{Proof of Theorem \ref{t5.5}} Let $f\in C[0,a]$. Then
$\{\langle Y_k^a(t),f\rangle\dvtx t\ge0\}$ converges a.s. to $\{
\langle Z
_t^a,f\rangle\dvtx t\ge
0\}$ in the topology of $D([0,\infty), \mbb{R})$. Consequently, we
have a.s. $\langle Y_k^a(t),f\rangle\to\langle Z_t^a,f\rangle$ for
a.e. $t\ge0$ [see,
e.g., Ethier and Kurtz (\citeyear{EtK86}), page 118]. By Lemma \ref{l5.1},
%
%
\begin{eqnarray}\label{5.8}
G(\langle Y_k^a(t),f\rangle)
&=&
G(\langle\lambda,f\rangle) + kb_k\int_0^t G^\prime(\langle Y
_k^a(s),f\rangle)
\langle\eta_k,f\rangle\,ds \nonumber\\
&&{}
- kb_k\int_0^t G^\prime(\langle Y_k^a(s),f\rangle)
\langle Y_k^a(s),f\rangle\,ds \nonumber\\
&&{}
+ \frac{1}{2}k^2\sigma_k^2\int_0^t
G^{\prime\prime}(\langle Y_k^a(s),f\rangle) \langle Y
_k^a(s),f^2\rangle\,ds \nonumber\\
&&{}
- \frac{1}{2}k\sigma_k^2\int_0^t G^{\prime\prime}(\langle Y
_k^a(s),f\rangle)
\langle Y_k^a(s),f\rangle^2 \,ds \\
&&{}
+ \int_0^tds\int_0^k m_k(dz)\int_{[0,a]} H(x,z,\langle Z
_s^a,f\rangle)
Y_k^a(s,dx) \nonumber\\
&&{}
+ \int_0^tds\int_0^k
[\varepsilon_k(s,z)+\xi_k(s,z)+\zeta_k(s,z)]m_k(dz) \nonumber\\
&&{}
+ \mbox{local mart.},\nonumber
\end{eqnarray}
where
\[
H(x,z,u) = G\bigl(u + zf(x)\bigr) - G(u) - G^\prime(u)zf(x)
\]
and
\[
\zeta_k(s,z) = \int_{[0,a]} [H(x,z,\langle Y
_k^a(s),f\rangle) -
H(x,z,\langle Z^a_s,f\rangle)] Y_k^a(s,dx).
\]
By the mean-value theorem,
\[
|\zeta_k(s,z)|\le\int_{[0,k]} |H^\prime_u(x,z,\theta_k(s))
\langle Y_k^a(s)-Z^a_s,f\rangle| Y_k^a(s,dx),
\]
where $\theta_k(s)$ takes values between $\langle Y_k^a(s),f\rangle$ and
$\langle Z^a_s,f\rangle$. For $G\in C^3(\mbb{R})$ we have
\begin{eqnarray*}
|H^\prime_u(x,z,\theta_k(s))|
&=&
\bigl|G^\prime\bigl(\theta_k(s) + zf(x)\bigr) - G^\prime(\theta_k(s)) -
G^{\prime\prime}(\theta_k(s))zf(x)\bigr| \\
&\le&
\|f\|\bigl(2\|G^{\prime\prime}\| + \tfrac{1}{2} \|f\|
\|G^{\prime\prime\prime}\|\bigr)(z\land z^2).
\end{eqnarray*}
It follows that
%
%
\begin{eqnarray}\label{5.9}
|\zeta_k(s)|
&\le&
\| f\|\bigl(2\|G^{\prime\prime}\| + \tfrac{1}{2} \|f\|
\|G^{\prime\prime\prime}\|\bigr)(z\land z^2) \nonumber\\[-8pt]\\[-8pt]
&&{}
\times\langle Y_k^a(s),1\rangle|\langle
Y_k^a(s)-Z^a_s,f\rangle|.\nonumber
\end{eqnarray}
By (\ref{5.9}) and Schwarz's inequality,
\begin{eqnarray*}
&&\mbf{E}\biggl[\int_0^{t\land\tau_n}ds\int_0^k |\zeta_k(s)|
m_k(dz)\biggr] \\
&&\qquad
\le C_k(t)\biggl\{\mbf{E}\biggl[\int_0^{t\land\tau_n}
\langle Y_k^a(s)-Z^a_s,f\rangle^2 \,ds\biggr]\biggr\}^{1/2} \\
&&\qquad\quad{}
\times\biggl\{\mbf{E}\biggl[\int_0^{t\land\tau_n}
\langle Y_k^a(s),1\rangle^2\,ds\biggr]\biggr\}^{1/2} \\
&&\qquad
\le\sqrt{n}C_k(t)\biggl\{\mbf{E}\biggl[\int_0^{t\land\tau_n}
\langle Y_k^a(s)-Z^a_s,f\rangle^2 \,ds\biggr]\biggr\}^{1/2},
\end{eqnarray*}
where
\[
C_k(t) =\|f\|\bigl(2\|G^{\prime\prime}\| + \tfrac{1}{2}
\|G^{\prime\prime\prime}\|\|f\|\bigr)\int_0^k (z\land z^2) m_k(dz).
\]
Note that $\sup_{k\ge1}C_k(t)< \infty$. It then follows that
\[
\mbf{E}\biggl[\int_0^{t\land\tau_n}ds\int_0^k |\zeta_k(s)|
m_k(dz)\biggr]\to0,\qquad k\to\infty.
\]
Now letting $k\to\infty$ in (\ref{5.8}) and using Lemmas \ref{l5.6}
and \ref{l5.7} we obtain (\ref{5.7}) for $G\in C^3(\mbb{R})$. A
simple approximation shows the martingale problem actually holds for
any $G\in C^2(\mbb{R})$.
\end{pf*}
\begin{theorem}\label{t5.8} Suppose that $kb_k\to b$, $\eta_k\to\eta$
weakly on $[0,a]$ and\break $k^2\sigma_k^2\delta_0(dz) + (z\land
z^2)m_k(dz)$ converges weakly on $[0,\infty)$ to a finite measure
$\sigma^2\delta_0(dz) + (z\land z^2)m(dz)$ as $k\to\infty$. Let
$\{0\le a_1< \cdots< a_n=a\}$ be an ordered set of constants. Then
$\{(Y_k^{a_1}(t), \ldots, Y_k^{a_n}(t))\dvtx t\ge0\}$ converges to
$\{(X_t^{a_1}, \ldots, X_t^{a_n})\dvtx t\ge0\}$ in distribution on
$D([0,\infty), M[0,a_1]\times\cdots\times M[0,a_n])$.
\end{theorem}
\begin{pf}
By Lemma \ref{l5.4} the sequence $\{(Y_k^{a_1}(t), \ldots,
Y_k^{a_n}(t))\dvtx t\ge0\}$ is tight in $D([0,\infty)$, $M[0,a_1]\times
\cdots\times M[0,a_n])$. Let $\{(Z_t^{a_1}, \ldots, Z_t^{a_n})\dvtx
t\ge0\}$ be a weak limit point of $\{(Y_k^{a_1}(t), \ldots,
Y_k^{a_n}(t))\dvtx t\ge0\}$. To get the result, we only need to show
$\{(Z_t^{a_1}, \ldots, Z_t^{a_n})\dvtx t\ge0\}$ and $\{(X_t^{a_1},
\ldots, X_t^{a_n})\dvtx t\ge0\}$ have identical distributions on
$D([0,\infty), M[0,a_1]\times\cdots\times M[0,a_n])$. By passing
to a subsequence and using Skorokhod's representation, we can assume
$\{(Y_k^{a_1}(t), \ldots,\allowbreak Y_k^{a_n}(t))\dvtx t\ge0\}$ converges to
$\{(Z_t^{a_1}, \ldots, Z_t^{a_n})\dvtx t\ge0\}$ almost surely in the
topology of $D([0,\infty)$, $M[0,a_1]\times\cdots\times
M[0,a_n])$. Theorem \ref{t5.5} implies $\{Z_t^{a_n}\dvtx\allowbreak t\ge0\}$ is an
immigration superprocess solving the martingale problem (\ref{5.7})
with $a=a_n$. Let $\bar{Z}_t^{a_i}$ denote the restriction of
$Z_t^{a_n}$ to $[0,a_i]$. Then $Z_t^{a_n} = \bar{Z}_t^{a_n}$ in
particular. We will show $\{(Z_t^{a_1}, \ldots, Z_t^{a_n})\dvtx t\ge
0\}$ and \mbox{$\{(\bar{Z}_t^{a_1}, \ldots, \bar{Z}_t^{a_n})\dvtx t\ge0\}$}
are indistinguishable. That will imply the desired result since
$\{(X_t^{a_1}, \ldots,\allowbreak X_t^{a_n})\dvtx t\ge0\}$ and
$\{(\bar{Z}_t^{a_1}, \ldots, \bar{Z}_t^{a_n})\dvtx t\ge0\}$ clearly
have identical distributions on $D([0,\infty), M[0,a_1]\times\cdots
\times M[0,a_n])$. By the general theory of c\`{a}dl\`{a}g processes,
the complement in $[0,\infty)$ of
\[
D(Z) := \{t\ge0\dvtx \mbf{P}(Z^{a_1}_t=Z^{a_1}_{t-}, \ldots,
Z^{a_n}_t=Z^{a_n}_{t-})=1\}
\]
is at most countable [see Ethier and Kurtz (\citeyear{EtK86}), page
131]. For any $t\in
D(Z)$ we have almost surely $\lim_{k\to\infty}Y_k^{a_i}(t) = Z_t^{a_i}$
for each $i=1,\ldots,n$ [see Ethier and Kurtz (\citeyear{EtK86}),
page 118]. By an
elementary property of weak convergence, for any $t\in D(Z)$ we almost
surely have
\begin{eqnarray*}
Z_t^{a_i}([0,a_i])
&=&
\lim_{k\to\infty}Y_k^{a_i}(t,[0,a_i])
=
\lim_{k\to\infty}Y_k^{a_n}(t,[0,a_i]) \\
&\le&
Z^{a_n}_t([0,a_i]) = \bar{Z}^{a_n}_t([0,a_i])
=
\bar{Z}^{a_i}_t([0,a_i]).
\end{eqnarray*}
Since Theorem \ref{t5.5} implies $\{Z_t^{a_i}\dvtx t\ge0\}$ is
equivalent to $\{\bar{Z}_t^{a_i}\dvtx t\ge0\}$, we have
\[
\mbf{E}[Z_t^{a_i}([0,a_i])]
=
\mbf{E}[\bar{Z}^{a_i}_t([0,a_i])].
\]
It then follows that almost surely
%
%
\begin{equation}\label{5.10}
\lim_{k\to\infty}Y_k^{a_i}(t,[0,a_i])
=
\bar{Z}^{a_i}_t([0,a_i]).
\end{equation}
On the other hand, since $Y_k^{a_n}(t)\to\bar{Z}_t^{a_n}$, for any
closed set $C\subset[0,a_i]$ we have
%
%
\begin{equation}\label{5.11}
\limsup_{k\to\infty}Y_k^{a_i}(t,C)
=
\lim_{k\to\infty}Y_k^{a_n}(t,C)
\le
\bar{Z}^{a_n}_t(C) = \bar{Z}^{a_i}_t(C).
\end{equation}
By (\ref{5.10}) and (\ref{5.11}) we have $Z^{a_i}_t = \lim_{k\to
\infty}Y_k^{a_i}(t) = \bar{Z}^{a_i}_t$. Then $\{Z_t^{a_i}\dvtx t\ge0\}$
and $\{\bar{Z}_t^{a_i}\dvtx t\ge0\}$ are indistinguishable since both
processes are c\`{a}dl\`{a}g.
\end{pf}

Let $\mcr{M}$ be the space of Radon measures on $[0,\infty)$ furnished
with a metric compatible with the vague convergence. The result of
Theorem \ref{t5.8} clearly implies the convergence of $\{Y_k(t)\dvtx
t\ge
0\}$ in distribution on $D([0,\infty),\mcr{M})$ with the Skorokhod
topology. {From} Theorem \ref{t5.8} we can also derive the following
generalization of a result of \citet{BeL06} [see also
Bertoin and Le Gall (\citeyear{BeL00}) for an earlier result].
\begin{corol}\label{c5.9} Suppose that $kb_k\to b$, $\eta_k\to\eta$
weakly on $[0,a]$ and $k^2\sigma_k^2\delta_0(dz) + (z\land
z^2)m_k(dz)$ converges weakly on $[0,\infty)$ to a finite measure
$\sigma^2\delta_0(dz) + (z\land z^2)m(dz)$ as $k\to\infty$. Let
$\{0\le a_1< \cdots< a_n\}$ be an ordered set of constants. Then
$\{(Y_k(t,a_1), \ldots, Y_k(t,a_n))\dvtx t\ge0\}$ converges to
$\{(X_t(a_1), \ldots, X_t(a_n))\dvtx t\ge0\}$ in distribution on
$D([0,\infty), \mbb{R}_+^n)$.
\end{corol}

\section*{Acknowledgment}

We are very grateful to the referee for  a careful reading of the paper
and helpful comments.

%

%
\printaddresses

\end{document}